\RequirePackage{fix-cm}
\documentclass[smallextended]{svjour3}       
\smartqed  

\usepackage{marvosym}
\usepackage{graphicx}
\begin{document}

\title{Sharp large deviation probabilities for sums of independent bounded random variables
}

\titlerunning{On Talagrand's inequalities and sharp large deviations}        

\author{Xiequan Fan   \and  Ion Grama  \and  Quansheng~Liu
}


\institute{X. Fan (\Letter) \and I. Grama  \and Q.~Liu \at
               Univ. Bretagne-Sud,  UMR 6205, LMBA, F-56000 Vannes, France \\
              \email{fanxiequan@hotmail.com; E-mail: ion.grama@univ-ubs.fr, quansheng.liu@univ-ubs.fr}       %
}
\date{Received: date / Accepted: date}

\maketitle

\begin{abstract}
We obtain some optimal inequalities on tail probabilities for sums of independent bounded random variables.
Our main result completes an upper bound on tail probabilities due to Talagrand by giving a one-term asymptotic expansion
for large deviations.  This result can also be regarded as sharp large deviations of types of Cram\'{e}r and Bahadur-Ranga Rao.


\end{abstract}

\section{Introduction}
Let $\xi _{1},...,\xi _{n}$ be a sequence of independent non-degenerate random variables (r.v.s) satisfying $\mathbf{E}\xi_{i}=0$.  Let $S_{n}=\sum_{i=1}^{n}\xi _{i}$. The study of sharp large deviation probabilities  has a long history. Many interesting asymptotic expansions have been established in Cram\'{e}r \cite{Cramer38}, Bahadur and Ranga Rao \cite{BR60}, Petrov \cite{Petrov75}  and Rozovky \cite{L05,L12}. Various exponential upper bounds have been obtained by Prohorov \cite{Pr59}, Nagaev \cite{N79}, Petrov \cite{Petrov95} and Talagrand \cite{Ta96}, see also McDiarmid \cite{M89}, Nagaev \cite{N03,N07} and \cite{F12,F13} for martingales. For a few results on lower bounds, we refer to Nagaev \cite{N02} and Rozovky \cite{L03}.

In this paper, we consider the sums of bounded from above r.v.s $ \xi_i \leq 1$. Let $\sigma _{i}^{2}=\mathbf{E}\xi _{i}^{2}$ and $\sigma^2=\sum_{i=1}^n\sigma_{i}^2$. In particular $\sigma= \sqrt{n}\, \sigma_{1}$ in the independent and identically
distributed (i.i.d.) case. The celebrated  Bennett inequality \cite{Be62} states  that:  if $|\xi_{i}|\leq 1$ for all $1\leq i \leq n$, then, for all $x>0$,
\begin{equation}\label{Beiq}
\mathbf{P}(S_{n}\geq x\sigma )\leq B(x,\sigma):=\left( \frac{x+\sigma }{\sigma }\right) ^{-\sigma
x-\sigma ^{2}}e^{x\sigma }.
\end{equation}
Bennett's inequality (\ref{Beiq}) is not optimal. One of the most well-known improvements on (\ref{Beiq}) is Hoeffding's inequality given by  (2.8) of \cite{Ho63},
which states that: if $\xi_{i} \leq 1$ for all $1\leq i \leq n$, then, for all $ x >0$,
\begin{eqnarray}
\mathbf{P}(S_{n}\geq x\sigma ) & \leq & H_{n}(x, \sigma) \label{Hoiq}\\
&:= & \left\{\left( \frac{\sigma }{
x+\sigma }\right)^{ x \sigma+\sigma^2 }\left( \frac{n}{n-x\sigma}\right)^{ n-x\sigma }
\right\}^{\frac{n}{n+\sigma^2} }\!\!\mathbf{1}_{\{x\leq \frac n \sigma\}}  \nonumber \\
&\leq&B(x,\sigma),\label{BeHiq}
\end{eqnarray}
where (and hereafter) by convention $(\frac{n}{n-x})^{n-x}=1$ when $x=n$.
By considering the following distribution
\begin{equation}\label{dis}
\mathbf{P}(\eta_{i}=1)=\frac{\sigma^2/n}{1+\sigma^2/n} \ \ \ \ \mbox{and}\ \ \ \  \mathbf{P}(\eta_{i}=-\sigma^2/n)=\frac{1}{1+\sigma^2/n},
\end{equation}
Hoeffding showed that (\ref{Hoiq}) is the best that can be obtained from the exponential Markov inequality
\[
\mathbf{P}( S_{n} \geq x\sigma ) \leq  \inf_{\lambda\geq0}\mathbf{E} e^{\lambda(S_{n}-x\sigma)},\ \  \ \ x \geq 0,
\]
since $\inf_{\lambda\geq0}\mathbf{E} \exp\{\lambda( \sum_i^n \eta_i-x\sigma)\}=H_{n}(x, \sigma)$ for all $0\leq x \leq \frac{n}{\sigma}$.

Hoeffding's inequality (\ref{Hoiq}) can be still improved. For sums of bounded random variables $(\xi_{i})_{i=1,...,n}$ satisfying $ -B\leq \xi_{i} \leq 1$ for some constant $B\geq 1$ and all $1\leq i \leq n$,
Talagrand \cite{Ta95} showed that, for all $0\leq x \leq c_1 \frac{\sigma}{ B}$,
\begin{eqnarray}
\mathbf{P}(S_{n}\geq x\sigma) &\leq& \left( \Theta(x)  + c_2 \frac{B}{\sigma}  \right)\inf_{\lambda\geq0}\mathbf{E}e^{\lambda(S_n-x\sigma)} \label{fska}\\
&\leq& \left( \Theta(x)  + c_2 \frac{B}{\sigma}  \right)H_{n}(x, \sigma),\label{Taie}
\end{eqnarray}
where $c_1, c_2 > 0$ are two absolute constants, $$\Theta(x)=\Big(1-\Phi(x) \Big) \exp\left\{\frac{x^2}{2} \right\}$$ is Mill's ratio up to a constant $\sqrt{2\pi}$ and $\Phi(x)=\frac{1}{\sqrt{2\pi}}\int_{-\infty}^{x}e^{-\frac {t^2}2}dt$ is the standard normal distribution function.
Since
\begin{eqnarray}
 \frac{1}{\sqrt{2\pi}(1+x)} \leq \Theta(x)\leq \frac{1}{\sqrt{\pi}(1+x)},\ \ \ \ \ x\geq 0, \label{mc}
\end{eqnarray}
(see \cite{F71}), Talagrand's inequality (\ref{Taie}) improves Hoeffding's inequality (\ref{Hoiq}) by adding a factor $\Theta(x)[1+o(1)]$ of order $ \frac1{1+x} $ in the range $0\leq x=o( \frac \sigma B)$ as $\frac B\sigma \rightarrow 0$.

The scope of this paper is to extend Talagrand's inequality (\ref{fska}) to a one-term asymptotic expansion similar to Cram\'{e}r \cite{Cramer38} and  Bahadur and Ranga Rao \cite{BR60}. We also gives some explicit expressions for the relation of $c_1$ and $c_2$ under a $(2+\delta)$th moment condition. In particular, Corollary \ref{co2.2} gives the following sharp large deviation result: if $\xi_i\leq 1$ and $ \mathbf{E}|\xi_{i}|^3\leq B \mathbf{E}\xi_{i}^2$ for some constant $B > 0$ and all $1\leq i \leq n$, then, for all $0 \leq x \leq 0.1 \frac{\sigma}{B}$,
\begin{eqnarray}\label{jkls}
\mathbf{P}(S_n\geq x\sigma ) &=& \left[  \Theta(x) +  16 \theta  \frac{B}{\sigma}  \right] \inf_{ \lambda \geq 0}\mathbf{E}e^{\lambda(S_n-x\sigma)},
\end{eqnarray}
where $|\theta|\leq 1$. Since $\theta \geq -1$, we complete Talagrand's upper bound (\ref{fska})  by giving a sharp lower bound.  This lower bound also improves the lower bound of Nagaev \cite{N02} for sums of bounded random variables.
Moreover, if $|\xi_{i}|\leq 1$, the constant $16$ in (\ref{jkls}) can be improved to $3.08$ (see Theorem \ref{th2.3}).
In the i.i.d. case, from (\ref{mc}) and (\ref{jkls}), we find that, for all $0 \leq x = o( \sqrt{n})$ as $n \rightarrow \infty$,
\begin{eqnarray}
\frac{\mathbf{P}(S_n \geq x\sigma )}{\Theta(x)\inf_{ \lambda \geq 0}\mathbf{E}e^{\lambda(S_n-x\sigma)}} = 1+ o(1).\label{efsa}
\end{eqnarray}

Notice that the expression $\inf_{ \lambda \geq 0}\mathbf{E}e^{\lambda(S_n-x\sigma)}$ in (\ref{jkls}) can be  rewritten in the form $\exp\{-n \Lambda_n^\ast(\frac{x\sigma}n)\}$, where $\Lambda_n^\ast(x)=\sup_{\lambda\geq 0}\{\lambda x - \frac1n \log \mathbf{E}e^{\lambda S_n}\}$ is the Fenchel-Legendre transform  of the cumulant function of $S_n$. The function $\Lambda^\ast(x)=\lim_{n\rightarrow \infty}\Lambda_n^\ast(x)$ is known as the good rate function in the large deviation principle (LDP) theory (see \cite{D98}).

To show the relation among equality (\ref{jkls}) and the results of Cram\'{e}r \cite{Cramer38} and  Bahadur and Ranga Rao \cite{BR60} when $n\rightarrow \infty$, we consider the i.i.d. case. In this case, equality (\ref{jkls}) reduces to the following results: for all $0\leq x \leq 0.1 \sqrt{n}\, \sigma_1/B$,
\begin{eqnarray}\label{cram1}
\mathbf{P}\left( \frac{S_{n}}{\sqrt{n} \sigma_1} \geq x \right)=   e^{- n \Lambda^\ast(x \sigma_1/\sqrt{n})} \Theta(x)\left[1 + O\left(\frac{1+x}{\sqrt{n}}\right) \right]
\end{eqnarray}
and, for all $0\leq y \leq 0.1 \sigma_1^2/B$,
\begin{eqnarray}\label{badrao1}
\mathbf{P}\left( \frac{S_{n}}{n} \geq y \right)=   e^{- n \Lambda^\ast(y/\sigma_1)}\left[ \Theta(y \sqrt{n}/\sigma_1) + O\left(\frac{1}{\sqrt{n}}\right) \right].
\end{eqnarray}
Cram\'{e}r \cite{Cramer38} showed the following equality:  for all $0\leq x =o (\sqrt{n})$,
\begin{eqnarray}\label{cram2}
\mathbf{P}\left( \frac{S_{n}}{\sqrt{n} \sigma_1} \geq x \right)= e^{- \frac12 x^2  + \frac{1}{\sqrt{n}} x^3 \lambda\left(\frac{x}{\sqrt{n}} \right)}   \Theta(x)\left[ 1+ O\left(\frac{1+x}{\sqrt{n}}\right) \right],
\end{eqnarray}
where $\lambda(\cdot)$ is the Cram\'{e}r series (see \cite{Petrov75} for details).
Bahadur and Ranga Rao \cite{BR60} obtained the following expansion: for given $y> 0$,
\begin{eqnarray}\label{badrao2}
\mathbf{P}\left( \frac{S_{n}}{n} \geq y \right)= \frac{ \sigma_1\mbox{}_y \, e^{- n \Lambda^\ast(y/\sigma_1)}}{t_y  \sqrt{2\pi n }  } \left[ 1+ O\left(\frac{c_1(y)}{n}\right) \right],
\end{eqnarray}
where $t_y$ and $\sigma_1\mbox{}_y$ depend on the distribution of $\xi_1$ and on $y$ in a complicated way. Compared to (\ref{cram2}) and (\ref{badrao2}), the equalities (\ref{cram1}) and (\ref{badrao1})
avoid the complicated factors $\lambda(\cdot)$, $t_y$ and $\sigma_1\mbox{}_y$. In particular, since
$$\Theta(x)=\frac{1}{x\sqrt{2\pi} }\Big[1- \frac{1}{x^2}+....  \Big]$$
 for $x\rightarrow \infty$, equality (\ref{badrao1}) implies that, for given  $y \in (0, 0.1 \sigma_1^2/B] $,
\begin{eqnarray}\label{badrao3}
\mathbf{P}\left( \frac{S_{n}}{n} \geq y \right)= \frac{ \sigma_1 \,  e^{- n \Lambda^\ast(y/\sigma_1)}}{y\sqrt{2\pi n }  } \left[ 1+ O\left(\frac{c_1(y)}{\sqrt{n}}\right) \right].
\end{eqnarray}
Compared to (\ref{badrao2}), equality (\ref{badrao3}) has the advantage that  the complicated factors $t_y$ and $\sigma_1\mbox{}_y$ have been replaced by the explicit values $y$ and $\sigma_1$, respectively.

In the proofs of our results, we make use of the conjugate distribution technique,
which becomes a standard for obtaining sharp large deviation expansions. This technique has been used  in Petrov \cite{Petrov75},  Nagaev \cite{N02}, Bercu and Rouault \cite{BR06},  Borovkov and Mogulskii \cite{BM07}, Petrov and Robinson \cite{P08}, Bercu, Coutin and Savyb \cite{Bs12} and  Gy\"{o}rfi,  Harrem\"{o}es and  Tusn\'{a}dy \cite{GHT12}.
Here we refine the technique inspired by Talagrand \cite{Ta95} and Grama and Haeusler \cite{GH00}.

The paper is organized as follows. In Section  \ref{sec2}, we present our main results.
In Section \ref{sec3}, we present some auxiliary results. In
Sections \ref{sec21} - \ref{sec23}, we prove the main results.

Throughout the paper, we make use of the following notation: $a\wedge b=\min\{a, b\}$, $a\vee b=\max\{a, b\}$, $a^+= a\vee0$, $\theta$ stands for a value satisfying $\left| \theta \right| \leq 1$.
We denote by $\mathcal{N}(0,1)$ the standard normal distribution, and we agree that $0^{-1}=\infty$.

\section{Main Results} \label{sec2}
Let $(\xi_{i})_{i=1,...,n}$ be a sequence of independent non-degenerate real random variables (r.v.s) with $\mathbf{E}\xi_{i}=0$ all over the paper.
In the sequel, we use the following condition.
\begin{description}
\item[(A)] There exist two constants $\delta \in (0,1]$ and $B>0$ such that, for all $\lambda\geq 0$,
\[
\mathbf{E}|\xi_{i}|^{2+\delta}\leq B^{2+\delta} \ \ \mbox{for all $1\leq i \leq n$ and} \ \ \ \ \  \sum_{i=1}^n \mathbf{E} \xi_{i}^2e^{\lambda \xi_{i}} \geq (1- B\lambda)\sum_{i=1}^n \mathbf{E} \xi_{i}^2.
\]
\end{description}

Notice that condition (A) is satisfied for $\delta=1$ if $\mathbf{E}|\xi_{i}|^3\leq B\mathbf{E}\xi_{i}^2$ holds. Indeed, since $(\mathbf{E}\xi_{i}^2)^{3/2}\leq  \mathbf{E}|\xi_{i}|^3 \leq B\mathbf{E}\xi_{i}^2$,
we have $\mathbf{E}\xi_{i}^2 \leq B^2$ and $\mathbf{E}|\xi_{i}|^3\leq B^3$. By the inequality $e^x\geq 1+x$, it follows that $\sum_{i=1}^n\mathbf{E} \xi_{i}^2e^{\lambda \xi_{i}} \geq \sum_{i=1}^n \mathbf{E} \xi_{i}^2(1+\lambda \xi_{i}) \geq (1- B\lambda) \sum_{i=1}^n \mathbf{E} \xi_{i}^2$ for all $\lambda\geq 0$.

For $ \delta \in (0, 1]$, denote by $C_{2+\delta}$ the Lyapunov constant defined as follows.  Suppose that $(\xi_i)_{i=1,...,n}$ have $(2+\delta)$th moments, i.e.  $\mathbf{E}|\xi_{i}|^{2+\delta}< \infty$ for all $1\leq i \leq n$. Then $C_{2+\delta}$ is the minimum of all absolute constants $C$ such that
 \[
  \sup_{x \in \mathbf{R}}\left| \frac{}{} \mathbf{P}(S_n \leq x\sigma )- \Phi(x)  \right|\leq C\, \frac{ \sum_{i=1}^{n} \mathbf{E} |\xi_{i}|^{2+\delta}}{\sigma^{2+\delta}}
 \]
holds for all $(\xi_{i})_{i=1,...,n}$. It is known that $0.4097 \leq  C_3 \leq 0.56$ and that $C_3 \leq 0.4784$ in the identically distributed case (see Shevtsova \cite{S10}). For the binomial distribution (for $0< p \leq 0.5$), Nagaev and Chebotarev \cite{N11} have recently proved that $C_3 \leq 0.4215$.

Our main result is the following theorem.
\begin{theorem}\label{th2.1}
 Assume  condition (A) and that $ \xi_{i} \leq1$ for all $1\leq i \leq n$. Then, for
 all $ 0\leq x < 0.25 \frac \sigma B$,
\begin{eqnarray}
  \mathbf{P}(S_n>x\sigma )  &=& \bigg(  \Theta(x) + \theta \varepsilon_x  \bigg)  \inf_{ \lambda \geq 0}\mathbf{E}e^{\lambda(S_n-x\sigma)}   \label{f17j} \\
  &\leq& \bigg( ( \Theta(x) +  \varepsilon_x  )\wedge 1 \bigg)\, H_n(x, \sigma),\label{f1sj}
 \end{eqnarray}
 where $|\theta|\leq1$ and
\begin{eqnarray*}
\varepsilon_x = \frac{e^{t}}{ 1- 2 t } \left(  \frac{1.58 }{\sqrt{\pi} } \frac{B}{\sigma} + \frac{2^{3+\delta} C_{2+\delta}  }{( 1- 2 t)^{\delta/2} } \frac{ \sum_{i=1}^{n} \mathbf{E} |\xi_{i}|^{2+\delta}}{\sigma^{2+\delta}} \right)
\end{eqnarray*}
with $t = \frac{2xB/\sigma}{1+\sqrt{ 1-4xB/\sigma }}$. In particular, in the i.i.d.\,case,  for all $ 0\leq x =o(n^{\delta/2})$,
\begin{eqnarray}
\frac{\mathbf{P}(S_n \geq x\sigma )}{\Theta(x)\inf_{ \lambda \geq 0}\mathbf{E}e^{\lambda(S_n-x\sigma)}} = 1+ o(1).\label{efsa}
\end{eqnarray}
and
\begin{eqnarray}
  \Big| \mathbf{P}(S_n>x\sigma ) - \Theta(x) \inf_{ \lambda \geq 0}\mathbf{E}e^{\lambda(S_n-x\sigma)} \Big| &=& O\left( \frac{1}{n^{\delta/2}}   H_n(x, \sigma)  \right)
\end{eqnarray}
as $n\rightarrow \infty$.
\end{theorem}

For r.v.s $\xi_i$ without moments of order larger than $2$, some  improvements of  Hoeffding's inequality (\ref{Hoiq})  can be found, for instance, in  Bentkus \cite{Be04} and Bentkus, Kalosha and van Zuijlen \cite{BZ06} and Pinelis \cite{P09}.  In particular,  when  $\xi_i \leq 1$ for all $1\leq i \leq n$,
Bentkus \cite{Be04} showed that
\begin{equation}\label{iebeq}
\mathbf{P}\left(S_{n}\geq x\right) \leq \frac{e^{2}}2\ \mathbf{P}^o \left(\sum_{i=1}^{n}\eta
_{i}\geq x \right),
\end{equation}
where $\eta _{i}$ are i.i.d. with distribution (\ref{dis}) and  $\mathbf{P}^o \left(\sum_{i=1}^{n}\eta
_{i}\geq x \right)$ is the log-concave hull of $\mathbf{P}(\sum_{i=1}^{n}\eta _{i}\geq x)$, i.e. $\mathbf{P}^o$ is the minimum log-concave function such that  $\mathbf{P}^o\geq \mathbf{P}$. As
$\eta_{i} \leq1$, inequality (\ref{iebeq}) is sharp up to an absolute constant $\frac{e^{2}}2$.
Here we give an equivalent to bound $\mathbf{P}^o \left(\sum_{i=1}^{n}\eta
_{i}\geq x \right)$. Applying (\ref{f17j}) to $\mathbf{P}^o \left(\sum_{i=1}^{n}\eta
_{i}\geq x \right)$ with $B=\max\{1, \frac{\sigma^2}n\}$, we find that if $\xi_i\leq 1$, then, for all $0\leq x \leq 0.24 \min \{\frac n \sigma, \sigma\}$,
\begin{eqnarray}
  \mathbf{P}^o \left(\sum_{i=1}^{n}\eta
_{i}\geq x \sigma \right) &=&  \Bigg(  \Theta(x) +  O\left(\max\left\{ \frac{1}{\sigma},  \frac{\sigma}{n} \right\} \right)   \Bigg) H_{n}(x,\, \sigma). \label{fsfgv}
\end{eqnarray}
Hence, we have the following inequality similar to (\ref{iebeq}) for $(\xi_{i})_{i=1,...,n}$ without moments of order larger than $2$.
\begin{corollary}\label{co2.1}
 Assume $\xi_{i}\leq1$ for all $1\leq i \leq n$. Then, for all $0\leq x \leq 0.24 \min \{\frac n \sigma, \sigma\}$,
\begin{eqnarray}
\mathbf{P}(S_n>x\sigma ) &\leq& F_1(x, \sigma)  H_n(x, \sigma), \label{dfcfv}
\end{eqnarray}
where $$F_1(x, \sigma)=\frac{e^2}{2}\Bigg(  \Theta(x) +  O\left(\max\left\{ \frac{1}{\sigma},  \frac{\sigma}{n} \right\} \right) \Bigg),\ \  \max\left\{ \frac{1}{\sigma},  \frac{\sigma}{n} \right\} \rightarrow 0.$$
\end{corollary}
This corollary shows us that (\ref{iebeq}) improves Hoeffding's bound $H_{n}(x,\, \sigma)$ by adding a factor $F_1(x, \sigma)$ in the range $0\leq x \leq 0.24 \min \{\frac n \sigma, \sigma\}$. Moreover, in the i.i.d. case, by (\ref{fsfgv}), we find that the ratio of bound (\ref{iebeq}) to bound (\ref{f1sj}) converges to $\frac{e^2}{2}$ for $0\leq x =o( n^{\delta/2})$ as $n\rightarrow \infty$, which means that (\ref{f1sj}) is better than (\ref{iebeq}) for all $0\leq x =o( n^{\delta/2})$. Of cause, the advantage of (\ref{dfcfv}) (also (\ref{iebeq}))  is that we do not assume that $\xi_i$ have moments of order larger than $2$ and that the missing factor exists in a larger range $0\leq x =o(\sqrt{n})$ as $n\rightarrow \infty$ in the i.i.d. case.

\vspace{0.1cm}

 Using Theorem \ref{th2.1} and $C_3 \leq 0.56$, we easily obtain  the following corollary.
\begin{corollary}\label{co2.2}
 Assume $\xi_{i}\leq1$ and $ \mathbf{E}|\xi_{i}|^3\leq B\mathbf{E}\xi_{i}^2$ for some constant $B>0$ and all $1\leq i \leq n$. Then, for all $0 \leq x \leq 0.1 \frac{\sigma}{B}$,
\begin{eqnarray}
\mathbf{P}(S_n>x\sigma ) &=& \left(  \Theta(x) +  16 \theta  \frac{B}{\sigma}  \right) \inf_{ \lambda \geq 0}\mathbf{E}e^{\lambda(S_n-x\sigma)}   \label{floer}\\
 &\leq& \Bigg(  \bigg( \Theta(x) + 16  \frac{B}{\sigma} \bigg)\wedge1 \Bigg)\, H_n(x, \sigma), \label{fhoe}
\end{eqnarray}
where $|\theta|\leq 1$.
\end{corollary}

Note that equality (\ref{floer}) implies Talagrand's inequality (\ref{fska}) by giving a large deviation expansion. In particular,
sice $\theta\geq -1$, equality (\ref{floer}) completes inequality (\ref{fska}) by giving a lower bound.

Some earlier lower bounds on tail probabilities, based on Cram\'{e}r large deviations, can be found in Arkhangelskii \cite{A89} and Nagaev \cite{N02}. In particular, Nagaev established the following lower bound
\begin{eqnarray}\label{cfbd}
\mathbf{P}(S_n>x\sigma ) &\geq& \Big( 1-\Phi(x) \Big) e^{-c_1x^3\frac{B}{\sigma} } \left( 1- c_2 (1+x)\frac{B}{\sigma} \right),
\end{eqnarray}
for some explicit constants $c_1, c_2$ and all $0 \leq x \leq \frac{1}{25} \frac{\sigma}{B}$. It is obvious that a precise lower bound of $\mathbf{E}e^{\lambda S_n }$ allows to improve Nagaev's bound  (\ref{cfbd}) by equality (\ref{floer}).

Inequality  (\ref{fhoe}) implies  the following Cram\'{e}r-type large deviations.
\begin{corollary}\label{co2.3}
 Assume $\xi_{i}\leq1$ and $ \mathbf{E}|\xi_{i}|^3\leq B\mathbf{E}\xi_{i}^2$ for some constant $B>0$ and all $1\leq i \leq n$. Then,
 for all $0 \leq x \leq 0.1 \frac{\sigma}{B}$,
\begin{eqnarray}
\mathbf{P}(S_n>x\sigma ) &\leq& \left(\frac{}{} 1-\Phi\left(  \check{x}
\right)\right)  \left[ 1+ 16 \sqrt{2 \pi}\left(1+  \check{x}
 \right) \frac{B}{\sigma} \right], \label{fbkl}
\end{eqnarray}
where $ \check{x} =\frac{x}{\sqrt{1+\frac{x}{3 \sigma} }}$ and satisfies
$$ \check{x} =x\left(1- \frac{x}{6\sigma}+o(\frac{x}{\sigma})\right)\ \ \ \ \ \  \mbox{as}  \ \ \ \ \ \ \ \ \frac x \sigma \rightarrow 0.$$
\end{corollary}

The interesting feature of the bound (\ref{fbkl}) is that it closely recovers the shape of the standard normal tail
for all $0\leq x = o( \sigma )$ as $ \sigma \rightarrow \infty$. Contrary to the Berry-Essen bound, the bound (\ref{fbkl})
is decreasing in an exponential rate for large $x$.

The well-known asymptotic expansions of tail probabilities (see Petrov \cite{Petrov75}) show that in general the value $ \check{x}$ in Corollary \ref{co2.3} can not be replaced by $x$. However, in the following sub-Gaussian case, $ \check{x}$ can be replaced by $x$.
\begin{theorem}\label{th2.2}
 Assume  $\xi_{i}\leq \sigma_{i}$ and $\mathbf{E}|\xi_{i}|^3\leq B\mathbf{E}\xi_{i}^2$ for some constant $B> 0$ and all $1\leq i \leq n$. Then, for
 all $ 0\leq x < 0.25 \frac \sigma B$,
\begin{eqnarray} \label{fco3}
\mathbf{P}(S_n>x\sigma ) &\leq&  \left(\frac{}{} 1- \Phi(x) \right) \left[1+ c_x \,(1+x ) \frac{B}{ \sigma } \right],
\end{eqnarray}
where
\[
c_x  = \frac{ e^{t+t^2} }{ 1- t }   \left(  \sqrt{2}    + \frac{ 16 \sqrt{2 \pi} C_3 e^{ \frac{1}{2}t^2 }  }{( 1- t )^{\frac 1 2}}\right)
\]
with $t = \frac{2xB/\sigma}{1+\sqrt{ 1-4xB/\sigma }}$.  In particular, in the i.i.d. case, for all $0 \leq x =o(\sqrt{n})$,
\begin{eqnarray}
\mathbf{P}(S_n>x\sigma ) &\leq&  \left(\frac{}{} 1- \Phi(x) \right) \Big[1+ o(1) \Big].
\end{eqnarray}
\end{theorem}

Note that the condition $\xi_{i}\leq \sigma_{i}$ (which replaces the condition $\xi_{i}\leq 1$ in Corollary \ref{co2.3}) is satisfied for
Rademacher r.v.s (i.e., $\mathbf{P}(\xi_{i}= \pm 1)=\frac12$).

For two-sided bounded r.v.s $|\xi_i|\leq B$ with $B>0$,   the following theorem shows that the constant $16$ in Corollary \ref{co2.2} can be further improved to a smaller one. Without loss of generality, we take $B=1$, otherwise we consider
$\xi_i/B$ instead of $\xi_i$.
\begin{theorem}\label{th2.3}
 Assume $|\xi_{i}|\leq 1$ for all $1\leq i \leq n$. Then,
 for all $0\leq x \leq 0.606\sigma$,
\begin{eqnarray}\label{sxcvy}
  \mathbf{P}(S_n>x\sigma )  &=& \left( \Theta(x) +   \theta \frac{ c_{x}}{\sigma}\right)  \inf_{ \lambda \geq 0}\mathbf{E}e^{\lambda(S_n-x\sigma)},
 \end{eqnarray}
where $|\theta|\leq 1$ and
\[
c_x=  \frac{2.24 e^{\frac{t^2}{2}} }{ \sqrt{ 1-t }} + \frac{ e^{t+t^2}}{\sqrt{\pi}(1-t) }
\]
with $t=  \frac{x}{\sigma} \exp\left\{ \frac{e x^2}{2 \sigma^2 } \right\}$.
In particular, for all $x\geq0$,
\begin{eqnarray}
  \mathbf{P}(S_n>x\sigma )
  &\leq& \bigg( \left(  \Theta(x) +  \frac{ c_{x}}{\sigma}  \right)\wedge 1  \bigg)\,H_n(x, \sigma).  \label{fhks}
 \end{eqnarray}
 Moreover, if $0\leq x \leq 0.1 \sigma$, then $c_x \leq 3.08$.
\end{theorem}

It is clear that inequality (\ref{fhks}) improves Hoeffding's  bound $H_n(x, \sigma)$ by adding a missing factor $\left(\Theta(x)+ c_x/\sigma \right)\wedge1$.

\begin{figure}
\includegraphics[width=0.48\textwidth]{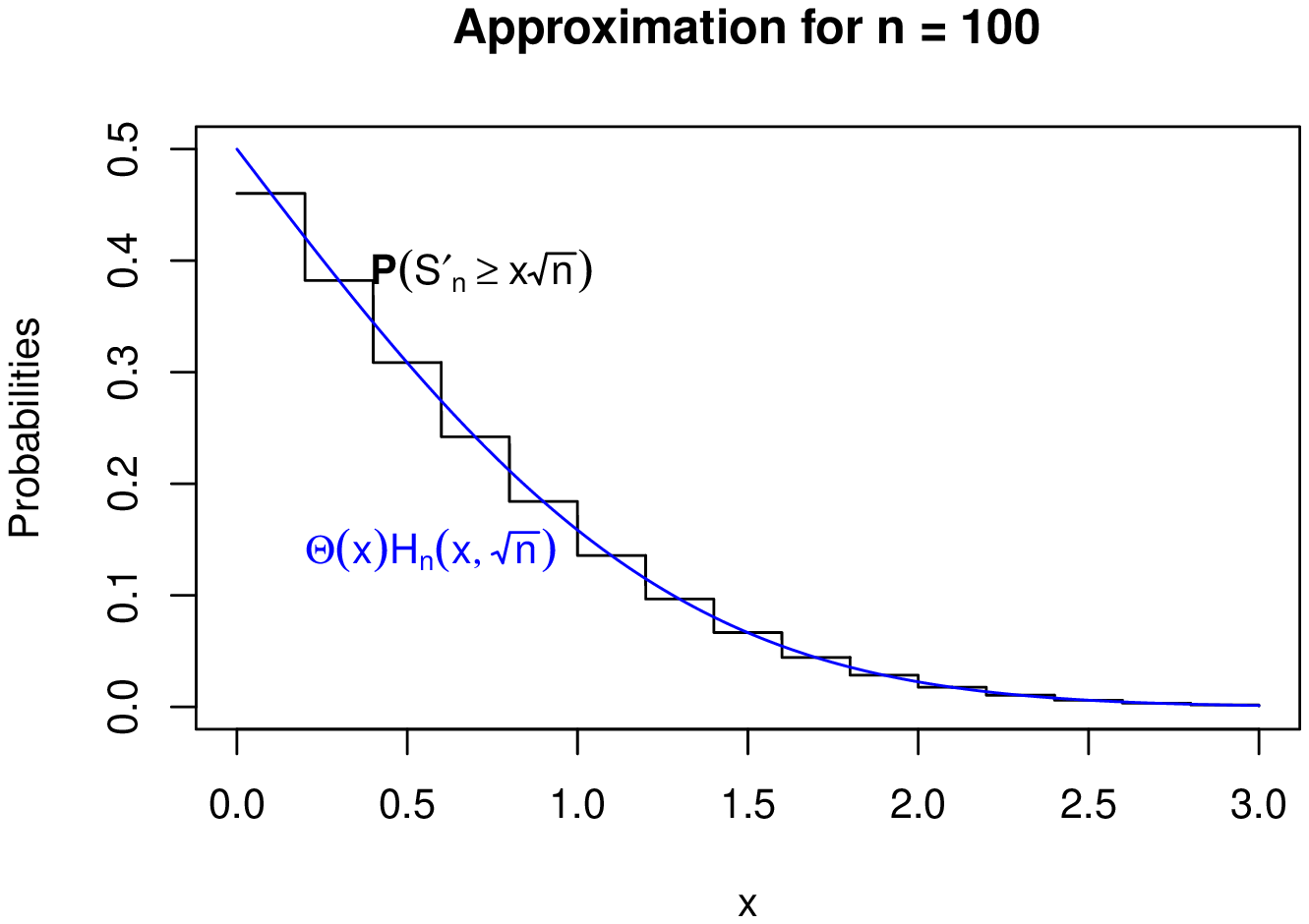}
\includegraphics[width=0.48\textwidth]{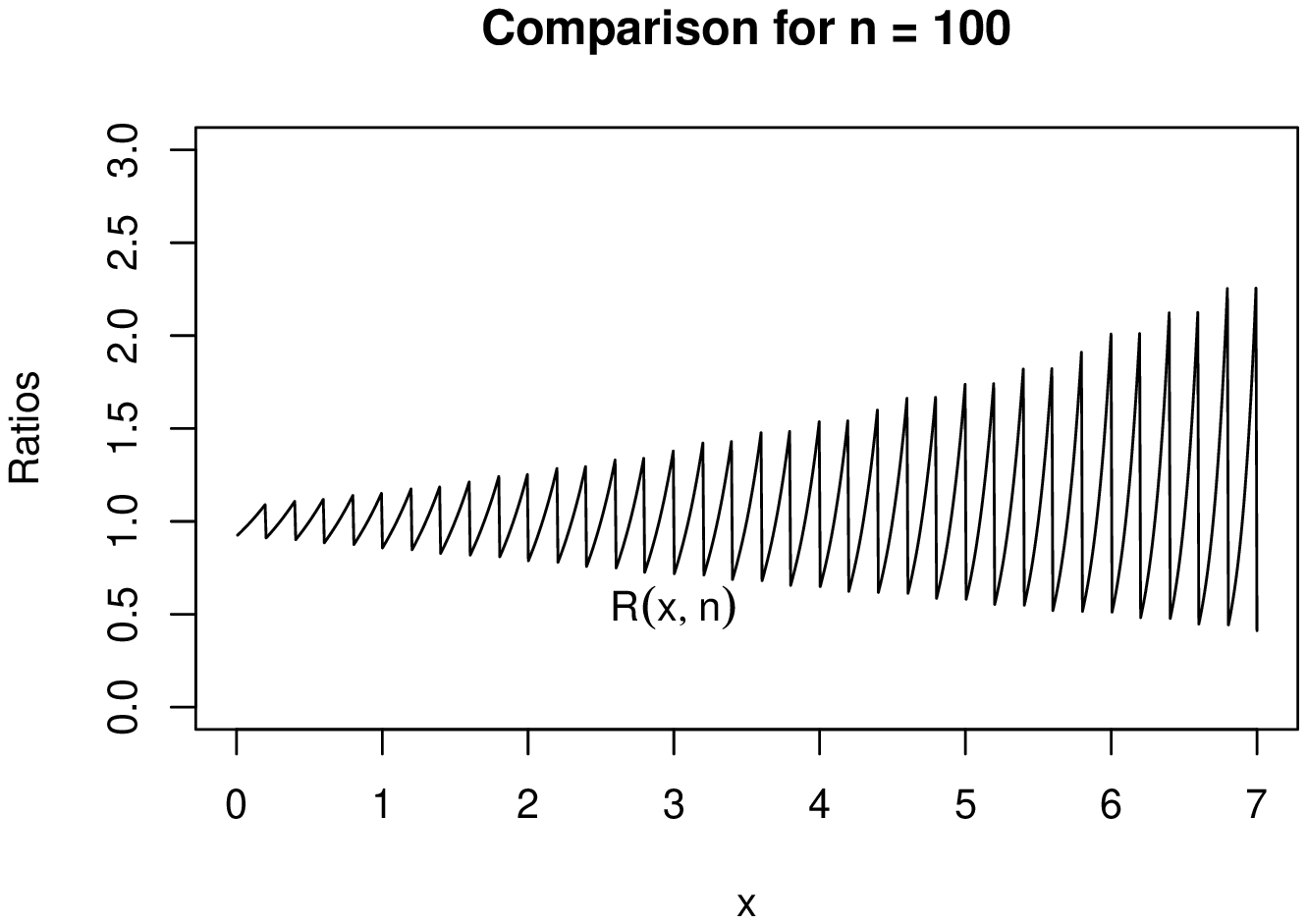}
\includegraphics[width=0.48\textwidth]{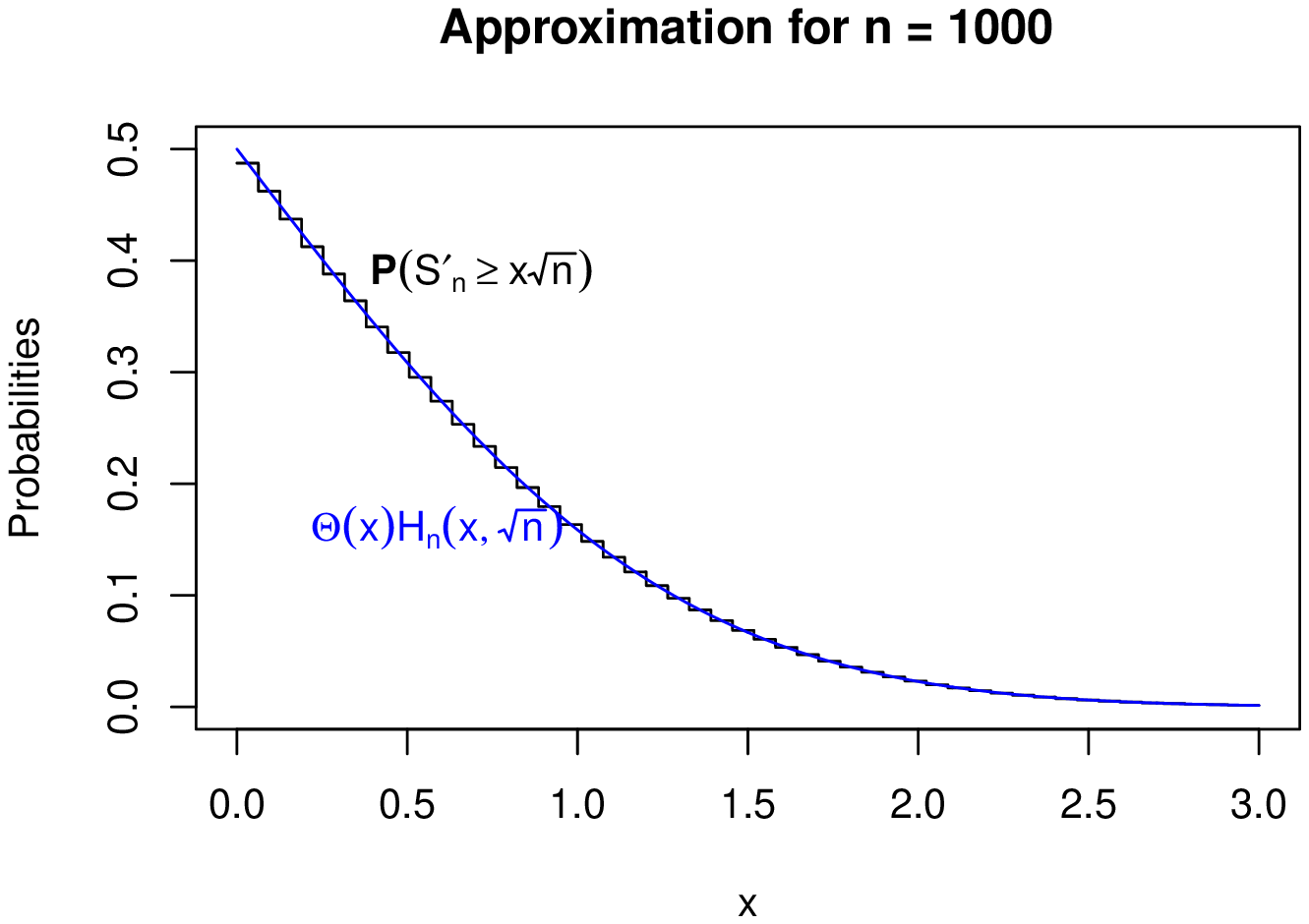}
\includegraphics[width=0.48\textwidth]{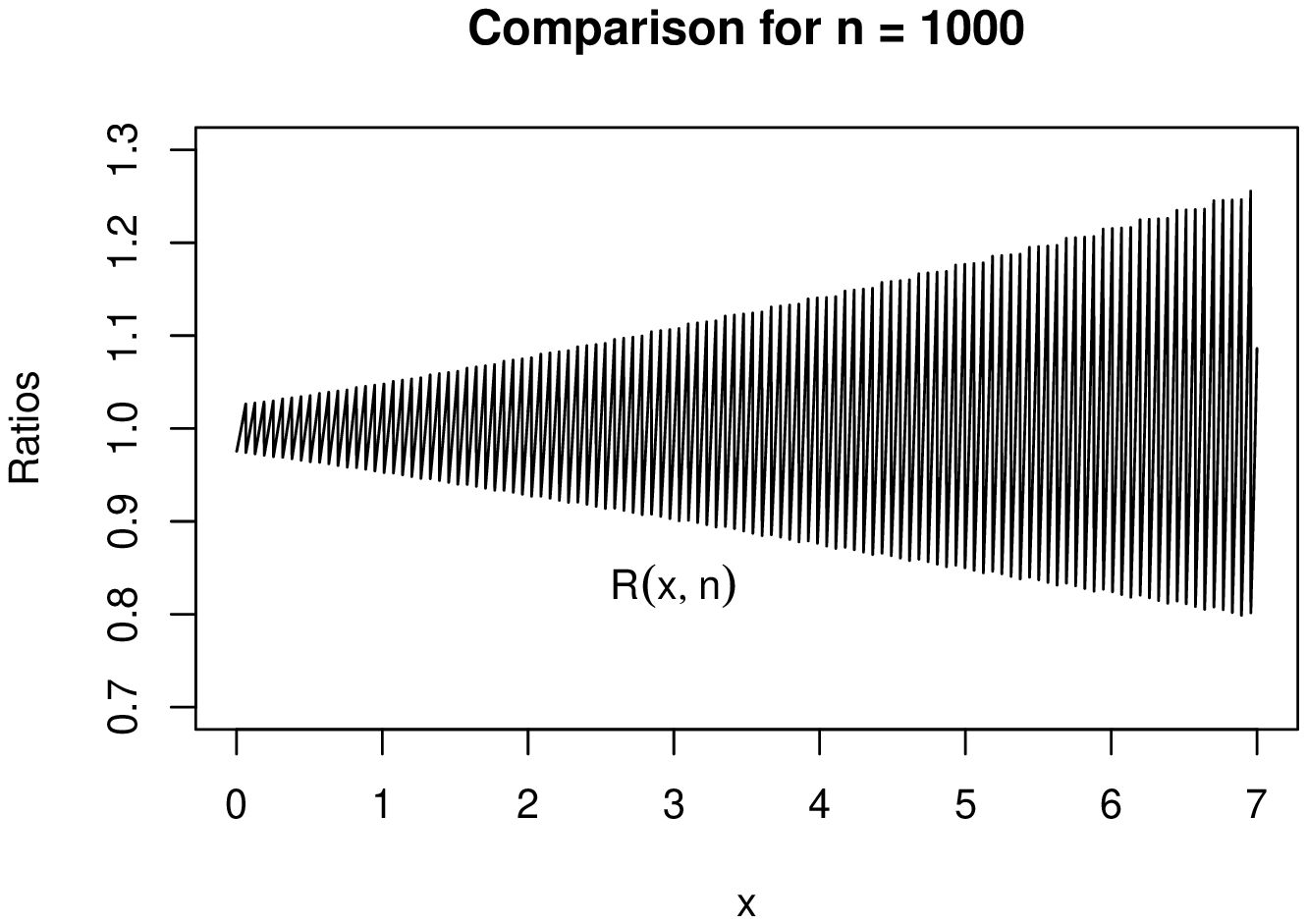}
\caption[]{Tail probabilities and ratios $R(x, n)$ for tail probabilities larger than $10^{-12}$ are displayed as a function of $x$ and various $n$.   }
\label{ssf}
\end{figure}

To show the tightness of equality (\ref{sxcvy}), let $S'_{n}=\varepsilon_1+...+\varepsilon_n$ be the sums of independent Rademacher r.v.s, i.e. $\mathbf{P}(\varepsilon_i=\pm 1)=\frac 12$ for all $i=1,...,n$. We display the simulation of $$R(x, n)=\frac{\mathbf{P}(S'_n \geq x \sqrt{n})}{\Theta(x)\inf_{ \lambda \geq 0}\mathbf{E}e^{\lambda(S'_n-x\sigma)}}=\frac{\mathbf{P}(S'_n \geq x \sqrt{n})}{\Theta(x)H_n(x, \sqrt{n})}  $$
in  Figure \ref{ssf}, which shows that $R(x, n)$ is very close to $1$ for large $n$'s.

\section{\textbf{Auxiliary Results}\label{sec3}}
Assume that $\mathbf{E}e^{\lambda \xi _i} < \infty$ for some constant $\lambda>0$ and all $i$.
 We consider the positive random variable
\[
Z_n(\lambda )=\prod_{i=1}^n\frac{e^{\lambda \xi _i}}{\mathbf{E}e^{\lambda \xi _i} },\ \ \ \ \ \ \ \lambda\geq 0, \label{C-1}
\]
so that $\mathbf{E}Z_{n}(\lambda)=1$ (the Esscher transformation). Introduce the
 \emph{conjugate probability measure} $\mathbf{P}_\lambda $ defined by
\begin{equation}
d\mathbf{P}_\lambda =Z_n(\lambda )d\mathbf{P}.  \label{dp}
\end{equation}
Denote by $\mathbf{E}_{\lambda}$ the expectation  with respect to $\mathbf{P}_{\lambda}$.
Setting
\[
b_i(\lambda)=\mathbf{E}_{\lambda} \xi_i=\frac{\mathbf{E}\xi
_ie^{\lambda \xi _i}}{\mathbf{E}e^{\lambda \xi _i}},\quad i=1,...,n,
\]
and
\[
\eta_i(\lambda)=\xi_i - b_i(\lambda), \quad i=1,...,n,
\]
we obtain the following decomposition:
\begin{equation}\label{xb}
S_k=B_k(\lambda )+Y_k(\lambda ),\quad k=1,...,n,
\end{equation}
where
\[
B_k(\lambda )=\sum_{i=1}^kb_i(\lambda ) \ \ \ \ \ \ \ \ \mbox{and}\ \ \ \ \ \  \ \  Y_k(\lambda )=\sum_{i=1}^k\eta _i(\lambda ).
\]
In the proofs of Theorems \ref{th2.1}-\ref{th2.3}, we shall need a two-sided bound of $B_{n}(\lambda)$. To this end,
we need some technical lemmas.

For a random variable bounded from above, the following inequality is well-known.
\begin{lemma}
\label{lemma3.1}  Assume $ \xi_{i} \leq1$.
Then, for all $\lambda \geq0,$
\begin{eqnarray*}
\mathbf{E} e^{\lambda \xi _i}  &\leq& Be(\lambda, \sigma_i^2),
\end{eqnarray*}
where
\[
Be(\lambda, t)=
\frac{  t }{1+t  }%
\exp\{\lambda\} + \frac{ 1}{1+ t }\exp\left\{-\lambda
t  \right\}.
\]
\end{lemma}

A proof of the inequality can be found in \cite{Be62}. This inequality is sharp and attains to equality when
\[
\mathbf{P}(\xi_i=1)=\frac{\sigma_i^2}{1+\sigma_i^2} \ \ \ \ \mbox{and}\ \ \ \ \mathbf{P}(\xi_i=-\sigma_i^2)=\frac{1}{1+\sigma_i^2}.
\]

By Lemma \ref{lemma3.1}, we easily obtain the following estimation  of the moment generating function $\mathbf{E}e^{\lambda  \xi _i}$.
\begin{lemma}
\label{lemma3.2}  Assume $ \xi_{i} \leq B$ and $\mathbf{E}|\xi_i|^{2+\delta}\leq B^{2+\delta}$ for some constants $B, \delta>0$
and all $1\leq i \leq n$.
Then, $\sigma_{i}^2 \leq B^2$ and, for all $\lambda \geq 0$,
\begin{eqnarray}
 \mathbf{E}e^{\lambda  \xi _i}  \leq e^{ \frac{B^2 \lambda^2}{2}}.
\end{eqnarray}
\end{lemma}

\noindent \textbf{\emph{Proof.}}
Using Jensen's inequality and $\mathbf{E}|\xi_i|^{2+\delta}\leq B^{2+\delta}$, we deduce
\[
\sigma_{i}^2=\mathbf{E} \xi _i^2 \leq (\mathbf{E} |\xi _i|^{2+\delta})^{2/(2+\delta)}\leq B^2.
\]
Since $\frac{d}{dt}Be(\lambda, t)  \geq 0$ for all $\lambda, t\geq 0$,
the function $Be(\lambda, t)$ is increasing in $t \geq 0$ for all $\lambda\geq0$. Hence, by Lemma \ref{lemma3.1} and the fact that $\frac{\sigma_i}{B}\leq 1$,  for all $\lambda\geq 0$,
\begin{eqnarray}
 \mathbf{E}e^{\lambda  \xi _i} =  \mathbf{E}e^{\lambda B  \frac{\xi _i}{B}}   \leq  \frac{e^{B\lambda}+e^{-B\lambda}}{2} \leq e^{ \frac{B^2 \lambda^2}{2}}  .
\end{eqnarray}
This completes the proof of Lemma \ref{lemma3.2}. \hfill\qed

In the following lemma, we give a two-sided bound for $B_{n}(\lambda)$.

\begin{lemma}
\label{lemma3.3}  Assume $ \xi_{i} \leq B$ for some constant $B>0$ and all $1\leq i \leq n$.
Then, for all $\lambda \geq0,$
\begin{eqnarray*}
B_n(\lambda ) &\leq& \frac{e^{B \lambda} -1 }{B}  \sigma^2.
\end{eqnarray*}
If $(\xi_i)_{i=1,...,n}$ satisfies condition (A), then, for all $ \lambda \geq 0,$
\begin{eqnarray}
B_n(\lambda )  \geq \left( 1- \frac{B\lambda}{2}  \right)\lambda \sigma^2e^{-\frac{B^2\lambda^2}{2}} .  \nonumber
\end{eqnarray}
\end{lemma}

\noindent \textbf{\emph{Proof.}}
By Jensen's inequality, we have, for all $\lambda\geq0$, $\mathbf{E}e^{\lambda\xi_{i}}\geq e^{\lambda \mathbf{E}\xi_{i}} =1$.
Since $\mathbf{E}\xi_ie^{\lambda\xi_i}=\mathbf{E}\xi_i(e^{\lambda\xi_i}-1)\geq0$ for $\lambda\geq0$, and $\xi_i\leq B$, we obtain the upper bound as follows: for all $\lambda\geq 0$,
\begin{eqnarray*}
B_{n}(\lambda)&\leq& \sum_{i=1}^n \mathbf{E} \xi_{i}e^{\lambda \xi_{i}} =\sum_{i=1}^n \int_{0}^{\lambda}\mathbf{E} \xi_{i}^2e^{t \xi_{i}}\, dt \\ &\leq & \sum_{i=1}^n \int_{0}^{\lambda} \sigma_{i}^2 \, e^{B t}  dt \\
&=&\frac{e^{B \lambda} -1 }{B}\, \sigma^2 .
\end{eqnarray*}
If $(\xi_{i})_{i=1,...,n}$ satisfies condition (A), it follows that, for all $\lambda\geq0$,
\begin{eqnarray}  \label{f27}
\sum_{i=1}^{n}\mathbf{E}\xi_{i} e^{\lambda \xi_{i}} &=& \int_{0}^{\lambda}\sum_{i=1}^{n} \mathbf{E}\xi_{i}^2 e^{t\xi_{i}}dt\\
&\geq &  \int_{0}^{\lambda}(1-  B t)dt \sum_{i=1}^{n}\mathbf{E}\xi_{i}^2 \nonumber \\
&=& \left( 1- \frac{B\lambda }{2} \right)\lambda\sigma^2.  \nonumber
\end{eqnarray}
Therefore, using Lemma \ref{lemma3.2}, we get the lower
bound of $B_n(\lambda )$: for all $  \lambda \geq 0,$
\begin{eqnarray}  \label{f28}
B_n(\lambda )= \sum_{i=1}^{n} \frac{\mathbf{E}\xi_{i} e^{\lambda \xi_{i}}}{\mathbf{E}e^{\lambda \xi_{i}}} \geq
 \left( 1- \frac{B\lambda}{2} \right )\lambda  \sigma^2e^{-\frac{B^2\lambda^2}{2}}, \nonumber
\end{eqnarray}
which completes the proof of Lemma  \ref{lemma3.3}.\hfill\qed

Next, we give an upper bound for the \emph{cumulant} function
\begin{equation}
\Psi _n(\lambda )=\sum_{i=1}^n\log \mathbf{E} e^{\lambda \xi _i}, \ \ \ \ \ \ \ \    \lambda \geq 0.
\label{C-3}
\end{equation}

\begin{lemma}
\label{lemma3.4}  Assume  $ \xi_{i} \leq1$ for all $1\leq i \leq n$.
Then, for all $\lambda \geq0,$
\begin{eqnarray*}
\Psi _n(\lambda ) &\leq& n \log\left( \frac{1}{1+\sigma^2/n}%
\exp\left\{-\lambda \sigma^2/n \right\} + \frac{\sigma^2/n}{1+\sigma^2/n} \exp\{\lambda\}\right).
\end{eqnarray*}
\end{lemma}

\noindent \textbf{\emph{Proof.}}
Since the function
\begin{eqnarray}
f(\lambda,t)= \log\left( \frac{1}{1+t}\exp\left\{-\lambda t \right\} +
\frac{t}{1+t} \exp\{\lambda\} \right),\ \ \ \ \ \lambda, t\geq0,  \nonumber
\end{eqnarray}
has a negative second derivative in $t>0$ (see Lemma 3 in  \cite{Ho63}),
then, for any fixed $\lambda\geq0$, $-f(\lambda,t)$ is convex in $t\geq 0$ and
\begin{equation}  \label{flt}
f(\lambda,t)\leq f(\lambda,0)+\frac{\partial}{\partial t}f(\lambda,0)\, t= (e^{\lambda}-1-\lambda)\, t,\ \ \ \ \
t\geq0.
\end{equation}
Therefore by Lemma \ref{lemma3.1} and Jensen's inequality, we get, for all $\lambda\geq 0$,
\begin{eqnarray*}
\Psi _n(\lambda ) &\leq & \sum_{i=1}^n f(\lambda, \sigma_{i}^2) \\
&\leq& n f(\lambda, \sigma^2/n) \\
 &=& n \log\left( \frac{1}{1+\sigma^2/n}%
\exp\left\{-\lambda \sigma^2/n \right\} + \frac{\sigma^2/n}{1+\sigma^2/n} \exp\{\lambda\}\right).
\end{eqnarray*}
This completes the proof of Lemma \ref{lemma3.4}.
\hfill\qed

Denote the  variance of $Y_n(\lambda )$  by $\overline{\sigma}^2(\lambda) = \mathbf{E}_{\lambda}Y_n^2(\lambda),\ \lambda\geq 0$. By the relation between $\mathbf{E}$ and $\mathbf{E}_{\lambda}$, the following inequality is obvious:
\begin{eqnarray*}
 \overline{\sigma}^2(\lambda) &= &\sum_{i=1}^n\left( \frac{\mathbf{E} \xi
_i^2e^{\lambda \xi _i}  }{\mathbf{E} e^{\lambda \xi _i} }-\frac{(\mathbf{E}\xi _ie^{\lambda \xi _i} )^2}{(\mathbf{E}e^{\lambda
\xi _i} )^2}\right)\nonumber,\ \ \ \lambda \geq0.
\end{eqnarray*}
The following lemma gives some estimations of $\overline{\sigma}^2(\lambda)$.
\begin{lemma}
\label{lemma3.5}
 Assume $ \xi_{i} \leq B$ and condition (A) for all $1\leq i \leq n$. Then,  for all $\lambda \geq0,$
\begin{eqnarray}
 (1- 2B\lambda)^+\sigma^2 \ \leq \ \overline{\sigma}^2(\lambda) \ \leq \  e^{B\lambda}\sigma^2.
\end{eqnarray}
\end{lemma}

\noindent \textbf{\emph{Proof.}}  Since $\mathbf{E}e^{\lambda\xi_{i}}\geq 1,  \lambda\geq 0,$ and $\xi_{i}\leq B$, we get, for all $\lambda\geq0,$
\begin{eqnarray*}
\overline{\sigma}^2(\lambda) \leq \sum_{i=1}^n \mathbf{E} \xi
_i^2e^{\lambda \xi _i} \leq \sum_{i=1}^n \mathbf{E} \xi_i^2e^{\lambda B} =  e^{B\lambda}\sigma^2.
\end{eqnarray*}
This gives the upper bound of $\overline{\sigma}^2(\lambda)$.
For all $\lambda\geq0$, it is easy to see that
\begin{eqnarray}
 \mathbf{E}\xi _ie^{\lambda \xi _i} = \int_{0}^\lambda \mathbf{E}\xi^2_ie^{t \xi _i}dt  \leq  \int_{0}^\lambda  e^{B t} \mathbf{E}\xi_{i}^2dt   =\left(\frac{e^{B\lambda}-1}{B} \right)\mathbf{E}\xi_{i}^2.
\end{eqnarray}
Using Lemma \ref{lemma3.2} and condition (A), we obtain, for all $  \lambda \geq0 $,
\begin{eqnarray*}
\overline{\sigma}^2(\lambda) &\geq&\sum_{i=1}^n \frac{\mathbf{E} \xi _i^2e^{\lambda \xi _i}   - (\mathbf{E}\xi _i e^{\lambda \xi _i} )^2}{e^{ B^2 \lambda^2}}\\
&\geq&   \frac{(1-B\lambda )\sigma^2 -  (e^{B\lambda} -1)^2 B^{-2}\sum_{i=1}^n (\mathbf{E}\xi_{i}^2)^2  }{e^{ B^2 \lambda^2} } \\
&\geq& e^{-B^2\lambda^2}\left( 1-  B\lambda  - (e^{B\lambda} -1)^2 \right)  \sigma^2 .
\end{eqnarray*}
Noting that $\overline{\sigma}^2(\lambda)\geq 0$, by a simple calculation, we have,
 for all $ \lambda \geq 0$,
\begin{eqnarray*}
\overline{\sigma}^2(\lambda) &\geq& (1- 2B\lambda)^+ \sigma^2,
\end{eqnarray*}
which gives the lower bound of $\overline{\sigma}^2(\lambda)$.
\hfill\qed

For the random variable $Y_n(\lambda )$, $\lambda\geq 0$, we have the following result on the rate of convergence to the standard normal law.
\begin{lemma}
\label{lemma3.6} Assume $ \xi_{i} \leq B$ and $\mathbf{E}|\xi_{i}|^{2+\delta} < \infty$ for some constants $ \delta \in(0, 1]$, $B>0$ and all $1\leq i \leq n$. Then,  for all $\lambda \geq0,$
\[
\sup_{y\in \mathbf{R}}\left| \mathbf{P}_\lambda \left(\frac{Y_n(\lambda )}{\overline{\sigma}(\lambda)} \leq y \right)-\Phi (y)\right| \leq \frac{2^{2+\delta} C_{2+\delta} e^{B\lambda} }
{\overline{\sigma}^{2+\delta}(\lambda)}\sum_{i=1}^{n} \mathbf{E} |\xi_{i}|^{2+\delta}.
\]
\end{lemma}

\noindent \textbf{\emph{Proof.}}   Notice that $Y_{n}(\lambda)=\sum_{i=1}^{n}\eta_{i}(\lambda)$ is the sum of independent r.v.s $\eta_{i}(\lambda)$
and $\mathbf{E}_{\lambda}\eta_{i}(\lambda)=0$. Using the well-known rate of convergence in the central limit theorem
(cf. e.g. \cite{Petrov75}, p. 115), we get, for all $\lambda \geq0$,
\begin{eqnarray*}
 \sup_{y\in \mathbf{R}}\left| \mathbf{P}_\lambda \left(\frac{Y_n(\lambda )}{\overline{\sigma}(\lambda)} \leq y\right )-\Phi (y)\right| \leq \frac{C_{2+\delta}}
{\overline{\sigma}^{2+\delta}(\lambda)}\sum_{i=1}^{n} \mathbf{E}_{\lambda} |\eta_{i}|^{2+\delta}.
\end{eqnarray*}
Using the inequality $(a+b)^{1+q} \leq 2^{q}(a^{1+q}+b^{1+q})$ for $a, b, q \geq0$,
we deduce, for all $\lambda\geq0$,
\begin{eqnarray*}
\sum_{i=1}^{n}\mathbf{E}_{\lambda} |\eta_{i}|^{2+\delta}  &\leq & 2^{1+\delta}\sum_{i=1}^{n} \mathbf{E}_{\lambda}(|\xi_{i}|^{2+\delta} +  |\mathbf{E}_{\lambda}\xi_{i}|^{2+\delta}  )\\
&\leq & 2^{2+\delta} \sum_{i=1}^{n} \mathbf{E}_{\lambda} |\xi_{i}|^{2+\delta}
\ \leq  \ 2^{2+\delta} \sum_{i=1}^{n} \mathbf{E} |\xi_{i}|^{2+\delta}e^{ \lambda\xi_{i} }  \\
&\leq & 2^{2+\delta} e^{B \lambda  }\,\sum_{i=1}^{n} \mathbf{E} |\xi_{i}|^{2+\delta} .
\end{eqnarray*}
Therefore, we obtain, for all $\lambda\geq0$,
\[
\sup_{y\in \mathbf{R}}\left| \mathbf{P}_\lambda \left(\frac{Y_n(\lambda )}{\overline{\sigma}(\lambda)} \leq y\right)-\Phi (y)\right| \leq  \frac{2^{2+\delta} C_{2+\delta} e^{B\lambda} }
{\overline{\sigma}^{2+\delta}(\lambda)}\sum_{i=1}^{n} \mathbf{E} |\xi_{i}|^{2+\delta}.
\]
This completes the proof of Lemma \ref{lemma3.6}.
\hfill\qed

We are now ready to prove the main technical result of this section.
\begin{theorem} \label{th3.1}
 Assume  $ \xi_{i} \leq B$ and $\mathbf{E}|\xi_{i}|^{2+\delta} < \infty$ for some constants $ \delta \in(0, 1]$, $B>0$ and all $1\leq i \leq n$. For an $x\geq 0$, if there exists a $\overline{\lambda}$ such that $\Psi _n'(\overline{\lambda} )=x\sigma$, then
\begin{eqnarray}
 \mathbf{P}(S_n>x\sigma )
 &=&  \bigg( \Theta\left(\overline{ \lambda} \overline{\sigma}(\overline{\lambda}) \right)  +  \theta \varepsilon_x   \bigg) \inf_{ \lambda \geq 0}\mathbf{E}e^{\lambda(S_n-x\sigma)},  \label{fmain}
\end{eqnarray}
where $|\theta|\leq 1$ and
\[
\varepsilon_x=  \frac{2^{3+\delta} C_{2+\delta}e^{B\overline{\lambda}}}{\overline{\sigma}^{2+\delta}(\overline{\lambda})} \sum_{i=1}^{n} \mathbf{E} |\xi_{i}|^{2+\delta}.
\]
\end{theorem}

\noindent \textbf{\emph{Proof}.}
According to the definition of the conjugate probability measure (cf. (\ref{dp})), we have the following representation of $\mathbf{P}(S_n>x\sigma )$: for given $ x, \lambda \geq0,$
\begin{eqnarray}\label{fff21}
\mathbf{P}(S_n>x\sigma )
&= & \mathbf{E}_\lambda (Z_n (\lambda)^{-1}\mathbf{1}_{\{S_n>x\sigma\}}) \nonumber\\
&= & \mathbf{E}_\lambda (e^{
-\lambda S_n+\Psi _n(\lambda ) }  \mathbf{1}_{\{S_n>x\sigma\}} )\nonumber\\
&= & \mathbf{E}_\lambda (e^{-\lambda x\sigma
 +\Psi _n(\lambda )-\lambda Y_n(\lambda)  -\lambda B_{n}(\lambda) +\lambda x\sigma  }   \mathbf{1}_{ \left\{Y_{n}(\lambda)+B_n(\lambda)-x\sigma >0 \right\} } )\nonumber\\
&= &e^{-\lambda x\sigma
 +\Psi _n(\lambda ) }  \mathbf{E}_\lambda (e^{-\lambda [Y_n(\lambda) + B_{n}(\lambda) -  x\sigma ] } \mathbf{1}_{ \left\{Y_{n}(\lambda)+B_n(\lambda)-x\sigma >0 \right\} }) \nonumber  .
\end{eqnarray}
Setting $U_{n}(\lambda)=\lambda (Y_n(\lambda) + B_{n}(\lambda) -  x\sigma ),$
we get
\begin{eqnarray} \label{fs26}
\mathbf{P}(S_n>x\sigma )   &=& e^{-\lambda x\sigma
 +\Psi _n(\lambda ) }   \int_0^\infty e^{- t}\mathbf{P}_{ \lambda  }(0< U_{n}(\lambda) \leq t) dt.
\end{eqnarray}
For an $x\geq 0$, if there exists a $\overline{\lambda}=\overline{\lambda}(x)$ such that $\Psi _n'(\overline{\lambda} )=x\sigma$, then
the exponential function $e^{-\lambda x\sigma +\Psi _n(\lambda ) }$ in (\ref{fs26}) attains its minimum at $
\lambda= \overline{\lambda}$.
Since $B_{n}(\overline{\lambda})=\Psi _n'(\overline{\lambda} )=x\sigma$, we have $U_{n}(\overline{\lambda})=\overline{\lambda} Y_n(\overline{\lambda})$ and
\begin{eqnarray} \label{jvb}
 e^{-\overline{\lambda} x\sigma
 +\Psi _n(\overline{\lambda} ) }
 &=& \inf_{\lambda\geq0} e^{ -\lambda x\sigma
 +\Psi _n(\lambda )  } =\inf_{ \lambda \geq 0}\mathbf{E}e^{\lambda(S_n-x\sigma)} .
\end{eqnarray}
Using Lemma \ref{lemma3.6}, we deduce
\begin{eqnarray}
 \int_0^\infty e^{- t}\mathbf{P}_{ \overline{\lambda}  }(0< U_{n}(\overline{\lambda}) \leq t) dt
&=&\int_0^\infty e^{- \overline{\lambda}   y \overline{\sigma}(\overline{\lambda})  }\mathbf{P}_{ \overline{\lambda}  } \left(0< U_{n}(\overline{\lambda}) \leq  \overline{\lambda}   y \overline{\sigma}(\overline{\lambda}) \right) \overline{\lambda}  \overline{\sigma}(\overline{\lambda})d y  \nonumber \\
&=&\int_0^\infty e^{- \overline{\lambda}   y \overline{\sigma}(\overline{\lambda})  }\mathbf{P} (0< \mathcal{N}(0,1) \leq   y  ) \overline{\lambda}  \overline{\sigma}(\overline{\lambda})d y  +  \theta \varepsilon_x \nonumber \\
&= &\int_0^\infty e^{ -  \overline{\lambda}  y \overline{\sigma}(\overline{\lambda}) }d\Phi\left(y \right) +  \theta \varepsilon_x  \nonumber \\
&=& \Theta\left(\overline{ \lambda} \overline{\sigma}(\overline{\lambda}) \right)  +   \theta \varepsilon_x, \label{f32}
\end{eqnarray}
where $|\theta|\leq 1$ and
\[
\varepsilon_x=\frac{ 2^{3+\delta} C_{2+\delta} e^{B\overline{\lambda}}}{\overline{\sigma}^{2+\delta}(\overline{\lambda})}  \sum_{i=1}^{n} \mathbf{E} |\xi_{i}|^{2+\delta} .
\]
Therefore, from (\ref{fs26}), for all $x\geq0$,
\begin{eqnarray}
 \mathbf{P}(S_n>x\sigma )
 =  \bigg( \Theta\left(\overline{ \lambda} \overline{\sigma}(\overline{\lambda}) \right)  +  \theta \varepsilon_x   \bigg) \inf_{ \lambda \geq 0}\mathbf{E}e^{\lambda(S_n-x\sigma)} . \nonumber  \end{eqnarray}
This completes the proof of Theorem \ref{th3.1}.\hfill\qed

\section{Proof of Theorem \ref{th2.1}} \label{sec21}
In the spirit of Talagrand \cite{Ta95}, we would like to make use of $\Theta(x)$ to approximate $\Theta(\overline{ \lambda} \overline{\sigma}(\overline{\lambda}))$ in Theorem \ref{th3.1}. The proof of Theorem \ref{th2.1}  is a continuation of the proof of Theorem \ref{th3.1}.

\noindent \textbf{\emph{Proof of Theorem \ref{th2.1}.}}
Using (\ref{jvb}) and Lemma \ref{lemma3.4}, we get, for all $x\geq0$,
\begin{eqnarray}
&& \inf_{ \lambda \geq 0}\mathbf{E}e^{\lambda(S_n-x\sigma)}  \nonumber\\
 &\leq& \inf_{\lambda \geq 0}  \exp\left\{-\lambda x\sigma
 +n \log\left( \frac{1}{1+\sigma^2/n}
\exp\left\{-\lambda \sigma^2/n \right\} + \frac{\sigma^2/n}{1+\sigma^2/n} \exp\{\lambda\}\right)  \right\}  \nonumber \\
&=& H_n(x,\sigma) . \label{kjls}
\end{eqnarray}
Since $|\Theta'(x)|\leq \frac{1}{\sqrt{\pi} x^2}$, we deduce
\begin{eqnarray}\label{fnsq}
\left|\Theta\left(\overline{ \lambda } \overline{\sigma}(\overline{\lambda}) \right) -\Theta(x) \right|&\leq & \frac{1}{\sqrt{\pi}}\frac{|x- \overline{ \lambda } \overline{\sigma}(\overline{\lambda})| }{ \overline{ \lambda }^2 \overline{\sigma}^2(\overline{\lambda}) \wedge x^2 } .
\end{eqnarray}
Using Lemma \ref{lemma3.3} and the inequality $e^{-x}\geq 1-x$ for  $x\geq0$, we have, for all $0\leq \overline{\lambda} \leq \frac{1}{B}$,
\begin{eqnarray}\label{ldkm}
 \left( 1-  B  \overline{\lambda} \right)\overline{\lambda}\sigma \leq \left( 1- \frac{B\overline{\lambda} }{2}\right)\overline{\lambda} e^{- \frac{B^2\overline{\lambda}^2}{2}} \sigma  \leq \frac{B_{n}(\overline{\lambda})}{\sigma }=x  \leq \frac{e^{B \overline{\lambda}} -1 }{B}  \sigma .
\end{eqnarray}
By the estimation of $\overline{\sigma}(\overline{\lambda})$ in Lemma \ref{lemma3.5}, it follows that, for all $0\leq \overline{\lambda} \leq \frac{1}{2B}$,
\begin{eqnarray}\label{fkssh}
 \left| x- \overline{ \lambda } \overline{\sigma}(\overline{\lambda}) \right|
&\leq& \overline{\lambda}\sigma \left [ \left( \frac{e^{B \overline{\lambda}} -1 }{B \overline{\lambda}}  - \sqrt{ 1- 2 B\overline{\lambda} }  \right) \vee \left( e^{\frac{B\overline{\lambda}}2}-(1-B\overline{\lambda})\right)  \right ]\nonumber \\
&\leq& { 1.58 e^{B \overline{\lambda}}}  B\overline{\lambda}^2 \sigma
\end{eqnarray}
and
\begin{eqnarray}
 \overline{ \lambda }^2 \overline{\sigma}^2(\overline{\lambda}) \wedge x^2
&\geq& \left( \overline{\lambda}^2\sigma^2  (1-2B\overline{\lambda}) \right) \wedge \left(  \overline{\lambda}^2\sigma^2 (1-B\overline{\lambda})^2 \right) \nonumber\\
&=&  \overline{\lambda}^2\sigma^2 (1-2B\overline{\lambda}) .
\end{eqnarray}
Hence, (\ref{fnsq}) implies that, for all $0\leq \overline{\lambda} < \frac{1}{2B}$,
\begin{eqnarray}\label{sg4f}
| \Theta\left(\overline{ \lambda } \overline{\sigma}(\overline{\lambda}) \right)-\Theta(x) | &\leq &  \frac{1.58 }{\sqrt{\pi}}\frac{e^{B\overline{\lambda}}}{(1-2B\overline{\lambda})  }  \frac{B}{\sigma}.
\end{eqnarray}
By  (\ref{ldkm}), it follows that $\left( 1-  B  \overline{\lambda} \right)\overline{\lambda} \leq \frac x  \sigma $ and
\begin{eqnarray}\label{landaba}
0\leq B\overline{\lambda} \leq t=: \frac{2xB/\sigma}{1+\sqrt{1-4Bx/\sigma}}\ \ \ \ \mbox{for all}\ \   \ 0\leq x \frac B \sigma < 0.25.
\end{eqnarray}
Therefore, by Lemma \ref{lemma3.5}, it is easy to see that
\begin{eqnarray}\label{ghjk}
 \frac{e^{B\overline{\lambda}}}{\overline{\sigma}^{2+\delta}(\overline{\lambda})}\, \sum_{i=1}^{n} \mathbf{E} |\xi_{i}|^{2+\delta}
 &\leq& \frac{e^{t} }{ ( 1-2  t) ^{1+ \frac {\delta}2} } \frac{ \sum_{i=1}^{n} \mathbf{E} |\xi_{i}|^{2+\delta}}{\sigma^{2+\delta}}.
\end{eqnarray}
Combining (\ref{fmain}) and (\ref{sg4f})-(\ref{ghjk}) together,
 we have, for all $ 0\leq  x < 0.25 \frac{\sigma}{B}$,
\begin{eqnarray*}
  \mathbf{P}(S_n>x\sigma ) &=&  \bigg( \Theta(x) + \theta \varepsilon_x \bigg)  \inf_{\lambda\geq 0}\mathbf{E}e^{\lambda(S_n-x\sigma)}\\
   &\leq& \bigg(\left( \Theta(x) +  \varepsilon_x \right)\wedge 1 \bigg) H_n(x,\sigma),
\end{eqnarray*}
where $|\theta|\leq 1$ and
\begin{eqnarray*}
\varepsilon_x = \frac{e^{t}}{ 1- 2 t } \left(  \frac{1.58 }{\sqrt{\pi} } \frac{B}{\sigma} + \frac{2^{3+\delta} C_{2+\delta}  }{( 1- 2 t)^{\delta/2} } \frac{ \sum_{i=1}^{n} \mathbf{E} |\xi_{i}|^{2+\delta}}{\sigma^{2+\delta}} \right).
\end{eqnarray*}
This completes the proof of Theorem \ref{th2.1}.\hfill\qed
\vspace{0.3cm}

\noindent \textbf{\emph{Proof of Corollary \ref{co2.3}.}} Since $\Theta(x)$ is decreasing in $x\geq0$, we deduce
$\Theta(x)\leq \Theta\left( \check{x}  \right)$ where $ \check{x} =\frac{x}{\sqrt{1+\frac{x}{3 \sigma} }}$.
Notice that Hoeffding's bound is less than Bernstein's bound, i.e.
\[
H_n(x,\sigma)\leq \exp\left\{ - \frac{ \check{x}^2}{2 } \right\}
\]
(cf. Remark 2.1 of \cite{F12}). Therefore, from (\ref{fhoe}), we have, for all $0\leq x \leq 0.1 \frac{\sigma}{B}$,
\begin{eqnarray}
  \mathbf{P}(S_n>x\sigma )  &\leq&  \left( \Theta\left( \check{x}  \right) + 16  \frac{B}{\sigma} \right)  \exp\left\{ -\frac{ \check{x} ^2}{2  }\right\} \nonumber \\
  &=& 1- \Phi \left( \check{x}  \right) +    16 \frac{B}{\sigma}   \exp\left\{ -\frac{ \check{x} ^2}{2  }\right\}.\nonumber
 \end{eqnarray}
 Using (\ref{mc}), we obtain, for all $0\leq x \leq 0.1 \frac{\sigma}{B}$,
 \begin{eqnarray*}
\mathbf{P}(S_n>x\sigma ) &\leq& \Big( 1-\Phi\left( \check{x}
\right)\Big)  \left[ 1+ 16 \sqrt{2 \pi} \left(1+  \check{x}
 \right) \frac{B}{\sigma} \right].
\end{eqnarray*}
This completes the proof of Corollary \ref{co2.3}.\hfill\qed

\section{Proof of Theorem \ref{th2.2}} \label{sec22}
Under the condition of Theorem \ref{th2.2}, we have the following upper bound of $\Psi_{n}(\lambda)$.
\begin{lemma}
\label{lemma5.1}
  Assume $\xi_{i}\leq \sigma_{i}$ for all $1\leq i \leq n$.
 Then, for all $\lambda \geq0,$
\begin{eqnarray*}
\Psi_n(\lambda )\leq  \frac{\lambda^2\sigma^2}{2}.
\end{eqnarray*}
\end{lemma}

\noindent \textbf{\emph{Proof.}}   Using Lemma \ref{lemma3.1}, we have, for all $t\geq0$,
\[
\mathbf{E} \exp\left\{t \frac{\xi_{i}}{\sigma_{i}} \right\}  \leq \frac{e^{-t}+e^t }{2}  \leq \exp\left\{ \frac{t^2}{2} \right\}.
\]
Hence, for all $\lambda\geq0$,
\begin{eqnarray*}
\mathbf{E} e^{\lambda \xi_{i}} \leq \exp\left\{ \frac{\lambda^2\sigma_{i}^2}{2}\right\} \ \ \ \ \ \mbox{and}\ \ \ \ \Psi_n(\lambda ) \leq  \frac{\lambda^2}{2}  \sum_{i=1}^n\sigma_{i}^2= \frac{\lambda^2\sigma^2}{2}.
\end{eqnarray*}
This gives the upper bound of $\Psi_n(\lambda )$.
\hfill\qed

Under the condition of Theorem \ref{th2.2}, we have the following  lower bound of $\overline{\sigma}^2(\lambda)$.
\begin{lemma}
\label{lemma5.2}
Assume $\xi_{i}\leq B$ and $\mathbf{E}|\xi_{i}|^3\leq B \mathbf{E}\xi_{i}^2$ for some constant $B> 0$ and all $1\leq i \leq n$.  Then,  for all $ \lambda \geq 0,$
\begin{eqnarray}
(1- B\lambda )^+ e^{-  B^2 \lambda^2} \sigma^2  \leq \overline{\sigma}^2(\lambda).
\end{eqnarray}
\end{lemma}

\noindent \textbf{\emph{Proof.}}
Denote by $f(\lambda)=\mathbf{E} \xi_i^2e^{\lambda\xi_{i}} \mathbf{E} e^{\lambda \xi _i} - (\mathbf{E} \xi _ie^{\lambda \xi _i})^2$.
It is easy to see that
\begin{eqnarray*}
f'(0) = \mathbf{E}\xi_i^3  \ \ \ \textrm{and} \ \ \
f''(\lambda)&=&\mathbf{E} \xi_i^4e^{\lambda\xi_{i}} \mathbf{E} e^{\lambda \xi _i} - (\mathbf{E}\xi _i^2e^{\lambda \xi _i})^2 \geq 0.
\end{eqnarray*}
Thus, for all $\lambda\geq 0$,
\begin{eqnarray}
f(\lambda) \geq f(0)+f'(0) \lambda = \mathbf{E}\xi_i^2+\lambda \mathbf{E}\xi_{i}^3.
\end{eqnarray}
Therefore,  for all $ \lambda \geq 0,$
\begin{eqnarray}
 \mathbf{E}_{\lambda}\eta_{i}^2  =  \frac{\mathbf{E} \xi _i^2e^{\lambda \xi _i} \mathbf{E} e^{\lambda \xi _i} -(\mathbf{E}\xi _i e^{\lambda \xi _i} )^2}{(\mathbf{E}e^{\lambda \xi _i} )^2}
 \geq  \frac{\mathbf{E}\xi_i^2+\lambda \mathbf{E}\xi_{i}^3}{(\mathbf{E}e^{\lambda \xi _i} )^2}  \nonumber .
\end{eqnarray}
Using $\mathbf{E}|\xi_{i}|^3\leq B \mathbf{E}\xi_{i}^2$,  we get, for all $ \lambda \geq 0,$
\begin{eqnarray}\label{fqg}
 \mathbf{E}_{\lambda}\eta_{i}^2  \geq  \frac{\mathbf{E} \xi_i^2+\lambda \mathbf{E} \xi_{i}^3}{(\mathbf{E}e^{\lambda \xi_i} )^2}
 \geq \frac{(1- B \lambda)\mathbf{E} \xi_i^2}{(\mathbf{E}e^{\lambda \xi_i} )^2} .
\end{eqnarray}
Since $\mathbf{E}|\xi_{i}|^3\leq B \mathbf{E}\xi_{i}^2$ implies that $\mathbf{E}|\xi_{i}|^3\leq B^3$, by Lemma \ref{lemma3.2} and $\overline{\sigma}^2(\lambda)\geq 0$, it follows that
\begin{eqnarray*}
\overline{\sigma}^2(\lambda) \geq \sum_{i=1}^n \frac{ (1- B\lambda  )^+ \mathbf{E} \xi_{i}^2}{ e^{  B^2 \lambda^2 }}  \geq  (1- B\lambda )^+ e^{-  B^2 \lambda^2} \sigma^2.
\end{eqnarray*}
This completes the proof of Lemma \ref{lemma5.2}.
\hfill\qed

Notice that the condition in Theorem \ref{th2.2} is stronger than that in Corollary \ref{co2.2}.
We can easily prove Theorem \ref{th2.2} by (\ref{f17j}) of Corollary \ref{co2.2}. However, in order to obtain a constant $c_x$ in Theorem \ref{th2.2} as small as possible, we make use of Theorem \ref{th3.1} to prove Theorem \ref{th2.2}.

\vspace{0.3cm}
\noindent \textbf{\emph{Proof of Theorem \ref{th2.2}.}}
By Lemma \ref{lemma5.1}, for all $x\geq0$, we have
\begin{eqnarray}
 \inf_{ \lambda \geq 0} \mathbf{E} e^{\lambda(S_n -x\sigma)} &=& \inf_{ \lambda \geq 0}\exp\left\{-  \lambda  x\sigma
 +\Psi _n(  \lambda   ) \right\}   \nonumber \\
  &\leq&  \inf_{ \lambda \geq 0} \exp\left\{-\lambda x\sigma
 + \frac{\lambda^2\sigma^2}{2} \right\} \nonumber \\
 &\leq&  \exp\left\{- \frac {x^2} 2 \right\}.
\end{eqnarray}
Since $\sigma_{i}^3\leq  \mathbf{E}|\xi_{i}|^3 \leq B\sigma_{i}^2$ and $\xi_{i}\leq \sigma_{i}$,
we have $\sigma_i \leq B$ and $\xi_{i}  \leq B$. Then, by Lemmas \ref{lemma3.5}, \ref{lemma5.2} and inequality (\ref{ldkm}), it follows that, for all $0\leq \overline{\lambda} < \frac{1}{B}$,
\begin{eqnarray}\label{fkh}
&&\left| x- \overline{ \lambda } \overline{\sigma}(\overline{\lambda}) \right|  \nonumber\\
&\leq& \overline{\lambda}\sigma \left [ \left( \frac{e^{B \overline{\lambda}} -1 }{B \overline{\lambda}}  - \sqrt{ 1- B\overline{\lambda}  }e^{-  \frac{B^2 \overline{\lambda}^2}{2}}   \right) \vee \left( e^{\frac{B\overline{\lambda}}2}- \left(1- \frac{B\overline{\lambda}}{2} \right)e^{-\frac{B^2\overline{\lambda}^2}{2}}\right)  \right ]\nonumber \\
&\leq& {  e^{B \overline{\lambda}}}  B\overline{\lambda}^2 \sigma
\end{eqnarray}
and
\begin{eqnarray}
 \overline{ \lambda }^2 \overline{\sigma}^2(\overline{\lambda}) \wedge x^2
&\geq& \left( \overline{\lambda}^2\sigma^2  (1- B\overline{\lambda})e^{-\frac{B^2 \overline{\lambda}^2 }{2}} \right) \wedge  \left( \overline{\lambda}^2\sigma^2 \left(1- \frac{B\overline{\lambda}}{2}\right)^2 e^{- B^2\overline{\lambda}^2 } \right)\nonumber\\
&\geq&  \overline{\lambda}^2\sigma^2 (1- B\overline{\lambda})^{1.3}.
\end{eqnarray}
Hence, by (\ref{fnsq}), we obtain, for all $0\leq\overline{\lambda}<\frac1B$,
\begin{eqnarray}\label{gjn}
| \Theta\left(\overline{ \lambda } \overline{\sigma}(\overline{\lambda}) \right)-\Theta(x) | &\leq & \frac{e^{B\overline{\lambda}}}{\sqrt{\pi}(1- B\overline{\lambda})^{1.3} }  \frac{B}{\sigma}.
\end{eqnarray}
By Lemma \ref{lemma5.2}, it is easy to see that, for all $0\leq\overline{\lambda}< \frac1B$,
\begin{eqnarray*}
 \frac{e^{B\overline{\lambda}}}{\overline{\sigma}^{3}(\overline{\lambda})}\, \sum_{i=1}^{n} \mathbf{E} |\xi_{i}|^{3}
 \leq  \frac{e^{ B\overline{\lambda}+\frac{3}{2}B^2\overline{\lambda}^2 }  }{( 1- B\overline{\lambda} )^{\frac 3 2}} \frac B \sigma .
\end{eqnarray*}
Thus, from Theorem \ref{th3.1},  it follows that, for all $0\leq x < 0.25 \frac \sigma  B$,
\begin{eqnarray*}
  \mathbf{P}(S_n>x\sigma )  &\leq& \left( \Theta(x)  + c_x  \frac{B}{\sigma } \right) \exp \left\{ -\frac{x^2}{2} \right\} \nonumber  \\
  &= &  1-\Phi(x)+ c_x  \frac{B}{\sigma } \exp \left\{ -\frac{x^2}{2} \right\},
\end{eqnarray*}
where
\[
 c_x   = \frac{ e^{t+t^2} }{ 1- t }   \left(  \sqrt{2}    + \frac{ 16 \sqrt{2 \pi} C_3 e^{ \frac{1}{2}t^2 }  }{( 1- t )^{\frac 1 2}}\right)
\]
with $t$ defined in (\ref{landaba}).
Since $\Theta(x)\geq \frac{1}{\sqrt{2\pi}(1+x)}, x\geq0$,
we have, for all $0\leq x < 0.25 \frac \sigma  B$,
\begin{eqnarray*}
  \frac{ \mathbf{P}(S_n>x\sigma ) }{1-\Phi(x)} &\leq&  1 + \sqrt{2\pi} c_x  (1+x) \frac{B}{\sigma }.
\end{eqnarray*}
In particular, if $0\leq x \leq 0.1 \frac \sigma  B $, using $C_3\leq 0.56$,  we have $\sqrt{2\pi}c_x \leq 32.47$.\hfill\qed

\section{Proof of Theorem \ref{th2.3}} \label{sec23}
For bounded r.v.s, the result on the convergence rate of $Y_{n}(\lambda)$ to the standard normal law (cf. Lemma \ref{lemma3.6}) can be improved to the following one.
\begin{lemma}
\label{lemma6.1} Assume $ |\xi_{i}| \leq 1$
for all $1\leq i \leq n$.
 Then,  for all $\lambda \geq0,$
\[
\sup_{y\in \mathbf{R}}\left| \mathbf{P}_\lambda \left(\frac{Y_n(\lambda )}{\overline{\sigma}(\lambda)} \leq y \right)-\Phi (y)\right| \leq  \frac{1.12 }
{\overline{\sigma} (\lambda)}.
\]
\end{lemma}

\noindent \textbf{\emph{Proof.}}   Using the well-known rate of convergence in the central limit theorem in \cite{Petrov75}, p. 115, and the fact that $C_{3}\leq 0.56$, we get, for all $\lambda \geq0$,
\begin{eqnarray*}
 \sup_{y\in \mathbf{R}}\left| \mathbf{P}_\lambda  \left(\frac{Y_n(\lambda )}{\overline{\sigma}(\lambda)} \leq y\right )-\Phi (y)\right| \leq  \frac{0.56 }{\overline{\sigma}^{3}(\lambda)}\sum_{i=1}^{n}\mathbf{E}_{\lambda} |\eta_{i}|^{3} .
\end{eqnarray*}
Since $|\eta_i|\leq |\xi_{i}|+\mathbf{E}_{\lambda}|\xi_{i}| \leq 2$, it follows that
$\sum_{i=1}^{n}\mathbf{E}_{\lambda} |\eta_{i}|^{3}  \leq  2 \overline{\sigma}^{2}(\lambda)$.
Therefore,  for all $\lambda\geq0$,
\[
\sup_{y\in \mathbf{R}}\left| \mathbf{P}_\lambda \left(\frac{Y_n(\lambda )}{\overline{\sigma}(\lambda)} \leq y\right)-\Phi (y)\right| \leq   \frac{1.12}
{\overline{\sigma} (\lambda)},
\]
which completes the proof of Lemma \ref{lemma6.1}.
\hfill\qed

\begin{lemma}
\label{lemma6.2}  Assume $ |\xi_{i}| \leq 1$ for all $1\leq i \leq n$.
Then, for all $\lambda \geq0,$
\begin{eqnarray}
B_n(\lambda )  \geq  (1-e^{-\lambda})e^{-\frac{\lambda^2}{2}}\sigma^2.  \nonumber
\end{eqnarray}
\end{lemma}

\noindent \textbf{\emph{Proof.}}
Since $\xi_{i}\geq -1$, for all $\lambda\geq0$,
\begin{eqnarray}
\sum_{i=1}^{n}\mathbf{E}\xi_{i} e^{\lambda \xi_{i}} &=& \sum_{i=1}^{n} \int_{0}^{\lambda}\mathbf{E}\xi_{i}^2 e^{t\xi_{i}}dt\\
&\geq &  \int_{0}^{\lambda}e^{-t}dt \sum_{i=1}^{n}\mathbf{E}\xi_{i}^2 \nonumber \\
&=& (1-e^{-\lambda})\sigma^2.  \nonumber
\end{eqnarray}
Using Lemma \ref{lemma3.2} with $B=1$, we get the lower
bound of $B_n(\lambda )$: for all $  \lambda\geq 0,$
\begin{eqnarray}
B_n(\lambda )= \sum_{i=1}^{n} \frac{\mathbf{E}\xi_{i} e^{\lambda \xi_{i}}}{\mathbf{E}e^{\lambda \xi_{i}}} \geq
 (1-e^{-\lambda})e^{-\frac{\lambda^2}2}\sigma^2. \nonumber
\end{eqnarray}
This completes the proof of Lemma \ref{lemma6.2}.
\hfill\qed

\noindent \textbf{\emph{Proof of Theorem \ref{th2.3}.}} The proof is close to the proof of Theorem \ref{th2.1}, but here $B=1$.
From (\ref{f32}), by Lemma \ref{lemma6.1}, we deduce
\begin{eqnarray}
&&\int_0^\infty e^{-t}\mathbf{P}_{\overline{\lambda}}(0<U_n(\overline{\lambda})\leq t )dt \nonumber\\
&=&\int_0^\infty e^{-\overline{\lambda} y\overline{\sigma}(\overline{\lambda})}\mathbf{P}_{\overline{\lambda}}\left(0<U_n(\overline{\lambda})\leq \overline{\lambda} y\sigma(\overline{\lambda}) \right)\overline{\lambda}  \sigma(\overline{\lambda})dy \nonumber \\
&=&\int_0^\infty e^{-\overline{\lambda} y\overline{\sigma}(\overline{\lambda})}\mathbf{P} \left(0<\mathcal{N}(0,1)\leq y \right)\overline{\lambda}  \sigma(\overline{\lambda})dy  + \theta\frac{2.24 }
{\overline{\sigma} (\overline{\lambda})} \nonumber \\
&= &\int_0^\infty e^{ -  \overline{\lambda}  y \overline{\sigma}(\overline{\lambda}) }d\Phi\left(y \right) + \theta\frac{2.24 }
{\overline{\sigma} (\overline{\lambda})} \nonumber \\
&=& \Theta\left(\overline{ \lambda} \overline{\sigma}(\overline{\lambda}) \right)  +  \theta\frac{2.24 }
{\overline{\sigma} (\overline{\lambda})},\nonumber
\end{eqnarray}
where $x\geq0$, $\overline{\lambda}$ satisfies $\Psi _n'(\overline{\lambda} )=B_n(\overline{\lambda})=x\sigma$  and $|\theta|\leq 1$.
Therefore, we obtain the following corresponding result of Theorem \ref{th3.1} for two-sided bounded r.v.s:
\begin{eqnarray}
\mathbf{P}(S_n>x\sigma )   &=&  \left( \Theta\left(\overline{ \lambda} \overline{\sigma}(\overline{\lambda}) \right)  + \theta\frac{2.24 }
{\overline{\sigma} (\overline{\lambda})}  \right)\inf_{ \lambda \geq 0}\mathbf{E}e^{\lambda(S_n-x\sigma)} \label{fhsd},
\end{eqnarray}
where $|\theta|\leq 1$.
Since $\overline{\lambda} e^{-\frac{\overline{\lambda}^2}2} \sigma^2  \leq (1-e^{-\overline{\lambda}})e^{- \frac{ \overline{\lambda}^2}2}\sigma^2 \leq B_{n}(\overline{\lambda})=x\sigma$ (cf. Lemma \ref{lemma6.2}), by Lemma \ref{lemma5.2}, we have, for all $0\leq \overline{\lambda} < 1$,
\begin{eqnarray}
 \overline{\lambda}\leq \frac x \sigma e^{ \frac{\overline{\lambda}^2}2} \leq  \frac x \sigma \exp\left\{ \frac{x^2e^{\overline{\lambda}^2 }}{2\sigma^2}  \right\}  \leq t=: \frac x \sigma \exp\left\{ \frac{e\, x^2}{2\sigma^2} \right\}   \label{sfdqfd}
\end{eqnarray}
and
\begin{eqnarray}\label{fghncs}
 \overline{ \lambda }^2 \overline{\sigma}^2(\overline{\lambda}) \wedge x^2
&\geq& \overline{\lambda}^2(1-\overline{\lambda}) e^{-\overline{\lambda}^2 }\sigma^2\wedge  \overline{\lambda}^2 e^{-\overline{\lambda}^2 }\sigma^2  \nonumber\\
&=& \overline{\lambda}^2(1-\overline{\lambda}) e^{-\overline{\lambda}^2 }\sigma^2 .
\end{eqnarray}
From (\ref{fnsq}), (\ref{fkh}) and (\ref{fghncs}), it follows that, for all $0\leq \overline{\lambda} < 1$,
\begin{eqnarray*}
\left| \Theta\left(\overline{ \lambda } \overline{\sigma}(\overline{\lambda}) \right) - \Theta(x) \right| &\leq& \frac{ e^{\overline{\lambda}+\overline{\lambda}^2}}{\sqrt{\pi}(1-\overline{\lambda}) } \frac{1}{\sigma} .
\end{eqnarray*}
Using  Lemma \ref{lemma5.2} again, we obtain, for all $0\leq \overline{\lambda} < 1$,
\begin{eqnarray}
 \frac{2.24 }
{\overline{\sigma} (\overline{\lambda})} \leq \frac{2.24 e^{ \frac {\overline{\lambda}^2}2} }{ \sqrt{ 1-\overline{\lambda}} } \frac{1}{\sigma}.
\end{eqnarray}
Returning to (\ref{fhsd}), we get, for all $0\leq \overline{\lambda} < 1$,
\begin{eqnarray}
\mathbf{P}(S_n>x\sigma )   &=& \left( \Theta\left(x\right)  + \theta  \frac{c_x(\overline{\lambda})}{\sigma}  \right)\inf_{ \lambda \geq 0}\mathbf{E}e^{\lambda(S_n-x\sigma)} \nonumber,
\end{eqnarray}
where $|\theta|\leq 1$ and
\[
c_x(\overline{\lambda})=  \frac{2.24 e^{ \frac {\overline{\lambda}^2}2} }{ \sqrt{ 1-\overline{\lambda}} }+ \frac{ e^{\overline{\lambda}+\overline{\lambda}^2}}{\sqrt{\pi}(1-\overline{\lambda}) }  .
\]
Noting that $c_{x}(\overline{\lambda})$ is increasing in $\overline{\lambda} \in [0, 1)$, we have, for all $0\leq x \leq 0.606\sigma$,
\begin{eqnarray}
\mathbf{P}(S_n>x\sigma )   &=& \left( \Theta\left(x\right)  + \theta  \frac{c_x}{\sigma} \right)\inf_{ \lambda \geq 0}\mathbf{E}e^{\lambda(S_n-x\sigma)} \nonumber,
\end{eqnarray}
where $|\theta|\leq 1$ and
\[
c_x=  \frac{2.24 e^{ \frac{ t^2}2} }{ \sqrt{ 1-t }} + \frac{ e^{t+t^2}}{\sqrt{\pi}(1-t) }
\]
with $t$ defined in (\ref{sfdqfd}).
In particular,  using the inequality $\mathbf{P}(S_n\geq x\sigma) \leq \inf_{ \lambda \geq 0}\mathbf{E}e^{\lambda(S_n-x\sigma)} \leq H_n(x, \sigma)$ (cf. (\ref{kjls})), we obtain, for all $x\geq0$,
\begin{eqnarray*}
\mathbf{P}(S_n>x\sigma ) &\leq& \bigg( \left(\Theta\left(x\right)  +  \frac{c_\alpha}{\sigma} \right)\wedge 1 \bigg) H_n(x, \sigma).
\end{eqnarray*}
This completes the proof of Theorem \ref{th2.3}.\hfill\qed



\begin{thebibliography}{9}
\bibitem{A89} Arkhangelskii, A.N. (1989) Lower bounds for probabilities of large deviations for sums of independent random variables. \emph{Theory Probab. Appl.} \textbf{34}, no 4, 565--575.


\bibitem{BR60} Bahadur, R. and Ranga Rao, R. (1960) On deviations of the sample mean. \emph{Ann. Math. Statist.} \textbf{31}, 1015--1027.

\bibitem{Be62} Bennett, G. (1962) Probability inequalities for sum of independent random variables.
\emph{J. Amer. Statist. Asso.} \textbf{57}, No. 297, 33--45.


\bibitem{Be04} Bentkus, V. (2004) On Hoeffding's inequality,
\emph{Ann. Probab.} \textbf{32}, No. 2, 1650--1673.

\bibitem{BZ06} Bentkus, V., Kalosha, N. and van Zuijlen, M. (2006) On domination of tail probabilities of (super)martingales: explicit bounds, \textit{Lithuanian. Math. J.} \textbf{46}, No. 1, 1--43.


\bibitem{BR06} Bercu, B. and Rouault, A. (2006) Sharp large deviations for the Ornstein-Uhlenbeck process,
\emph{Theory Probab. Appl.} \textbf{46}, No. 1, 1--19.


\bibitem{Bs12} Bercu, B., Coutin, L. and Savyb, N. (2012) Sharp large deviations for the non-stationary Ornstein-Uhlenbeck process, \emph{Stochastic Process. Appl.} \textbf{122}, 3393--3424.


\bibitem{BM07} Borovkov, A.A., Mogulskii, A.A. (2007) On large and superlarge deviations of sums of independent random vectors under Cramer's condition. II, \emph{Theory Probab. Appl.} \textbf{51}, 567--594.

\bibitem{Cramer38} Cram\'{e}r, H. (1938) Sur un nouveau th\'{e}or\`{e}me-limite de la th\'{e}orie des probabilit\'{e}s. \emph{%
Actualite's Sci. Indust.} \textbf{736}, 5--23.


\bibitem{D98} Dembo, A. and Zeitouni, O. (1998) \emph{Large deviations techniques and applications} (2nd ed.). Springer, New
York.


\bibitem{F12} Fan, X., Grama, I. and Liu, Q. (2012)
 Hoeffding's inequality for supermartingales. \emph{Stochastic Process. Appl.} \textbf{122} 3545--3559.

\bibitem{F13} Fan, X., Grama, I. and Liu, Q. (2013) Cram\'{e}r large deviation expansions for martingales under Bernstein's condition,
 \textit{Stochastic Process. Appl.} \textbf{123}, 3919--3942.

\bibitem{F71} Feller, W. (1971) \emph{An introduction to
probability theory and its applications.} J. Wiley and Sons.

\bibitem{GH00} Grama, I. and Haeusler, E. (2000) Large
deviations for martingales via Cramer's method.
\emph{Stochastic Process. Appl.} \textbf{85}, 279--293.

\bibitem{GHT12}  Gy\"{o}rfi, L., Harrem\"{o}es, P. and  Tusn\'{a}dy, G. (2012) Some refinements of large deviation tail probabilities.
\emph{arXiv:1205.1005v1 [math.ST]}.

\bibitem{Ho63} Hoeffding, W. (1963) Probability inequalities
for sums of bounded random variables. \emph{J. Amer. Statist. Assoc.} \textbf{58}, 13--30.


\bibitem{M89} McDiarmid, C. (1989) On the method of bounded differences, in Surveys in Combina-
torics, ed J. Siemons, London Mathematical Society Lecture Note Series 141, Cambridge
University Press, 1989.

\bibitem{N79} Nagaev, S.V. (1979) Large deviations of sums
of independent random variabels. \textit{Ann. Probab.} \textbf{7}, No. 5,
745--789.

\bibitem{N02} Nagaev, S.V. (2002) Lower bounds for the probabilities of large deviations of sums of independent random
variables. \emph{Theory Probab. Appl.} \textbf{46}, no 1, 79--102; no 4, 728--735.


\bibitem{N03} Nagaev, S.V. (2003) On probability and moment
inequalities for supermartingales and martingales. \textit{Acta. Appl. Math.}
\textbf{79}, 35--46.

\bibitem{N07} Nagaev, S.V. (2007) On probability and moment
inequalities for supermartingales and martingales. \textit{Acta. Appl. Math.}
\textbf{97}, 151--162.

\bibitem{N11} Nagaev, S.V. and Chebotarev, V.I. (2011) On an estimate for the closeness of the  binomial distribution to the the normal distribution. \emph{Dokl. Math.} \textbf{83}, no 1, 19--21.

\bibitem{L03} Rozovky, L.V. (2003) A lower bound of large-deviation probabilities for the sample mean under the Cram\'{e}r condition. \emph{J. math. Sci.} \textbf{118}, No 6.

\bibitem{L05} Rozovky, L.V. (2005) large deviation probabilities for some classes of distributions statisfying the Cram\'{e}r condition. \emph{J. math. Sci.} \textbf{128}, No 1.

\bibitem{L12} Rozovky, L.V. (2012) Superlarge deviation probabilities for sums of independent lattice random variables with exponential decreasing tails. \emph{Statist. Probab. Letter} \textbf{82}, 72--76.


\bibitem{Petrov75} Petrov, V.V. (1975) \emph{Sums of
Independent Random Variables.} Springer-Verlag. Berlin.


\bibitem{Petrov95} Petrov, V.V. (1995) \textit{Limit Theorems of Probability Theory.}
Oxford University Press, Oxford.


\bibitem{P08} Petrov, V.V. and Robinson, J., (2008) Large deviations for sums of independent non identically distributed random variables.
\textit{Comm. Statist. Theory Methods} 37, 2984--2990.

\bibitem{P09} Pinelis, I. (2009) On the Bennett-Hoeffding inequality, \textit{arXiv:0902.4058v1}.

\bibitem{Pr59} Prohorov, Yu.V. (1959) An extremal problem in
probability theory. \textit{Theor. Probability Appl.} \textbf{4}, 201--203.

\bibitem{S10} Shevtsova, I.G. (2010) An improvement of convergence rate estimates
in the Lyapunov theorem. \emph{Doklady. Math.} \textbf{82} 862--864.

\bibitem{Ta95} Talagrand, M. (1995) The missing factor in
Hoeffding's inequalities. \emph{Ann. Inst. H. Poincar\'{e}  Probab. Statist.}  \textbf{31},
689--702.
\bibitem{Ta96} Talagrand, M. (1996) A new look at independence. \emph{Ann. Probab.}  \textbf{22}, 1--34.

\end{thebibliography}
\end{document}